\documentclass[12pt, a4paper]{article}

\usepackage[utf8]{inputenc}
\usepackage[T1]{fontenc}
\usepackage[english]{babel}
\usepackage{indentfirst}
\usepackage{amsmath}
\usepackage{amsfonts}
\usepackage{amssymb}
\usepackage{amsthm}
\usepackage{graphicx}
\usepackage{babelbib}
\usepackage{lmodern} \normalfont
\usepackage{color}
\usepackage{dsfont}
\usepackage[displaymath, mathlines, pagewise]{lineno}

\usepackage{fancyhdr}
\newtheorem{thm}{Theorem}[section]
\newtheorem{prop}[thm]{Proposition}
\newtheorem{lem}[thm]{Lemma}
\newtheorem{cor}[thm]{Corollary}
\newtheorem{assump}[thm]{Assumption}

\def\p{{\bf p}}
\def\q{{\bf q}}

\def\P{{\bf P}}
\def\Q{{\bf Q}}
\def\hp{\hat{{\bf p}}}
\def\hq{\hat{{\bf q}}}
\def\hP{\hat{{\bf P}}}
\def\hQ{\hat{{\bf Q}}}
\def\hg{\hat{g}}
\def\hf{\hat{f}}
\def\hh{\hat{h}}
\newcommand{\IN}{\mathbb{N}}
\newcommand{\kP}{\mathop{\mathfrak{P}}}
\newcommand{\IP}{\mathbb{P}}
\newcommand{\IE}{\mathbb{E}}
\newcommand{\mE}{\mathfrak{E}}
\newcommand{\Z}{\mathbb{Z}}

\newcommand{\R}{\mathbb{R}}
\newcommand{\F}{{\mathfrak{F}}}
\newcommand{\A}{{\mathsf{A}}}
\newcommand{\B}{{\mathcal{B}}}
\newcommand{\BB}{{\mathsf{B}}}
\newcommand{\D}{{\mathsf{D}}}
\newcommand{\C}{{\mathsf{C}}}
\newcommand{\mfS}{{\mathfrak{S}}}
\newcommand{\I}{{\mathcal{I}}}
\newcommand{\II}{{\mathfrak{I}}}
\newcommand{\M}{{\mathcal{M}}}

\newcommand{\mS}{{\mathcal{S}}}
\newcommand{\V}{{\mathcal{V}}}
\newcommand{\W}{{\mathcal{W}}}

\newcommand{\mH}{{\mathfrak{H}}}

\newcommand{\ind}{\mathds{1}}
\newcommand{\bdelta}{\boldsymbol{\delta}}
\newcommand{\eps}{\varepsilon}
\newcommand{\0}{\textbf{0}}

\let\phi=\varphi

\newcommand{\eqlaw}{\stackrel{\text{\tiny law}}{=}}

\newcommand{\1}[1]{{\mathds 1}{\{#1\}}}

\title{On uniform closeness of local times of Markov chains
and i.i.d.\ sequences}

\date{}

\author{Diego F. de Bernardini \and Christophe Gallesco \and 
Serguei Popov \thanks{Department of Statistics, Institute of Mathematics, 
Statistics and Scientific Computation, University of Campinas --
UNICAMP, rua S\'ergio Buarque de Holanda 651,
13083--859, Campinas SP, Brazil;
e-mails: $\{$bernardini, gallesco, popov$\}$@ime.unicamp.br}}

\begin{document}
	
\maketitle

\begin{abstract}
\footnotesize 
In this paper we consider 
 the field of local times of
 a discrete-time Markov chain
on a general state space, and obtain uniform (in time) upper
bounds on the total variation distance between this field
and the one of a sequence of~$n$ i.i.d.\ 
random variables with law given by the invariant measure of 
that Markov chain. The proof of this result
uses a refinement of the soft local time method of~\cite{PT15}.

\vspace{0.3cm}

\end{abstract}

\section{Introduction}
\label{intro}
The purpose of this paper is to compare the field of
 local times of a discrete-time 
Markov process with the corresponding field of i.i.d.\ 
random variables distributed
according to the stationary measure of this process,
 in total variation distance.
Of course, local times (also called occupation times)
of Markov processes is a very well studied subject.
It is frequently possible to obtain a complete
characterization of the law of this field
in terms of some Gaussian random field or process,
especially in continuous time (and space) setup.
The reader is probably familiar with 
Ray-Knight theorems as well as Dynkin’s and Eisenbaum’s 
isomorphism theorems; cf.\ e.g.\ \cite{R,Szn12}.
One should observe, however, that these theorems
usually work in the case when the underlying Markov
process is reversible and/or symmetric in some sense,
something we do not require in this paper.

To explain what we are doing here,
let us start by considering the following example: 
let $(X_j)_{j\geq 1}$ be a Markov chain 
on the state space $\Sigma=\{0,1\}$, 
with the following transition probabilities:
$\IP[X_{n+1}=k\mid X_n=k]=1-\IP[X_{n+1}=1-k\mid X_n=k]
=\frac{1}{2}+\eps$ for $k=0,1$, 
where~$\eps\in(0,\frac{1}{2})$ is small. 
Clearly, by symmetry, $(\frac{1}{2},\frac{1}{2})$ is
the stationary distribution of this Markov chain. 
Next, let $(Y_j)_{j\geq 1}$
be a sequence of i.i.d.\ Bernoulli random variables 
with success probability~$\frac{1}{2}$.
What can we say about the distance in total variation between 
the laws 
of $(X_1,\ldots,X_n)$ and $(Y_1,\ldots,Y_n)$? Note that the 
``na\"\i{}ve'' way of trying to force
the trajectories to be equal (given $X_1=Y_1$, use the maximal
 coupling of $X_2$ and $Y_2$;
if it happened that $X_2=Y_2$, then try to couple $X_3$ 
and $Y_3$, and so on) works
only up to $n=O(\eps^{-1})$. Even though this method is probably
not optimal, in this case it is easy to obtain that the total
variation distance converges to~$1$ 
as $n\to\infty$. This is because of the
 following: consider the event 
\[
\Xi^Z = \Big\{\frac{1}{n}\sum_{j=1}^{n-1} 
\ind_{\{Z_j = Z_{j+1}\}} > \frac{1}{2} + \frac{\eps}{2} \Big\},
\]
where $Z$ is~$X$ or~$Y$.
Clearly, the random variables $\ind_{\{Z_j = Z_{j+1}\}}$, 
$j\in\{1,\ldots,n-1\}$ are i.i.d.\ Bernoulli,
with success probabilities $\frac{1}{2}+\eps$ and $\frac{1}{2}$ 
for $Z=X$ and~$Z=Y$
correspondingly. Therefore, if $n\gg \eps^{-2}$, 
it is elementary to obtain that that
$\IP[\Xi^X]\approx 1$ and $\IP[\Xi^Y]\approx 0$, 
and so the total variation distance
between the \emph{trajectories} of~$X$ and~$Y$ 
is almost~$1$ in this case.

So, even in the case when the Markov chain gets quite close 
to the stationary distribution 
just in one step, usually it is not possible to couple its trajectory 
with an i.i.d.\  sequence,
unless the length of the trajectory is relatively short.
Assume, however, that we are not interested in the exact 
trajectory of~$X$ or~$Y$,
but rather, say, in the number of visits to~$0$ up to time~$n$. 
That is, denote
\[
L_n^Z(0) = \sum_{j=1}^n \ind_{\{Z_j=0\}}
\]
for $Z=X$ or $Y$. Are $L_n^X(0)$ and $L_n^Y(0)$ 
close in total variation distance
for \emph{all}~$n$?

Well, the random variable $L_n^Y(0)$ has the binomial 
distribution with parameters~$n$
and~$\frac{1}{2}$, so it is approximately Normal with 
mean~$\frac{n}{2}$ 
and standard deviation~$\frac{\sqrt{n}}{2}$. As for $L_n^X(0)$, 
it is elementary to obtain that it is 
approximately Normal with mean~$\frac{n}{2}$ and
standard deviation~$\sqrt{n}\big(\frac{1}{2}+O(\eps)\big)$.
Then, it is also elementary to obtain that the total variation distance between
 these two Normals is~$O(\eps)$, \emph{uniformly} in~$n$
(indeed, that total variation distance equals the total variation distance
between the Standard Normal and the centered Normal with variance 
$(1+O(\eps))^2$; that distance is easily verified to be of order~$\eps$).
This \emph{suggests} that the total variation distance between~$L_n^X(0)$
and~$L_n^Y(0)$ should be also of order~$\eps$ uniformly in~$n$.
Observe, by the way, that the distribution of the local times
of a two-state Markov chain can be explicitly written 
(cf.~\cite{BG}), so one can obtain a rigorous proof 
of the last statement in a direct way, after some work.  

Let us define the \emph{local time} of a stochastic process~$Z$
at site~$x$ at time~$n$ as the number of visits
to~$x$ up to time~$n$:
\begin{linenomath}
\begin{equation*}
L^Z_n(x) = \sum_{j=1}^n \ind_{\{Z_j=x\}}
\end{equation*}
\end{linenomath}
(sometimes we omit the  upper index when it is clear which process 
we are considering).
The above example shows that, if one is only interested in the local 
times of the
Markov chain (and not the complete trajectory),
 then there is hope to obtain a coupling
with the local times of an i.i.d.\ random sequence 
(which is much easier to handle).
Observe that there are many quantities of interest that 
can be expressed in terms 
of local times only
(and do not depend on the order), such as, for instance,
\begin{itemize}
\item hitting time of a site~$x$: $\tau(x) = \min\{n: L_n(x)>0\}$;
\item cover time: $\min\{n: L_n(x)>0 \text{ for all }x\in\Sigma\}$, 
where $\Sigma$ is the space where the process lives;
\item blanket time~\cite{DLP}:  $\min\{n\geq 1: L_n(x)\geq \delta n \pi(x)\}$,
where $\pi$ is the stationary measure of the process 
and $\delta\in (0,1)$ is a parameter;
\item disconnection time~\cite{DSzn,Szn10}:
loosely speaking, it is the time~$n$ when the set $\{x: L_n(x)>0\}$ 
becomes ``big enough'' to ``disconnect''
the space~$\Sigma$ in some precise sense;
\item the set of favorite (most visited) sites (e.g.~\cite{HS,T}): 
$\{x: L_n(x)\geq L_n(y)\text{ for all }y\in\Sigma\}$;
\item and so on.  
\end{itemize}
This justifies the importance of finding couplings as above.
Note also that, although 
not every Markov chain comes close to 
the stationary distribution in just one step, that 
can be sometimes circumvented by considering 
the process at times $k, 2k, 3k, \ldots$ with a large~$k$.  

In particular, we expect that our results would be useful 
when dealing with
\emph{excursion processes} (i.e., when~$\Sigma$ is a set
 of excursions of a random walk).
One may e.g.~refer to \cite{MS}, cf.\ Lemma~2.2 there
(note that the order of excursions does not matter,
so one would be able to get rid of the factor~$m$
in the right-hand side); 
also, we are working now on applications of our results to 
the decoupling for random
interlacements~\cite{BGP}.
\begin{figure}
 \centering \includegraphics{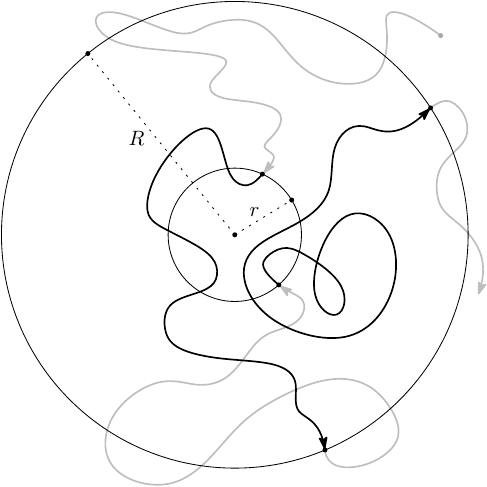} 
\caption{Excursions between the smaller and the 
larger spheres are pictured as bold pieces of the trajectory.}
\label{f_excursions}
\end{figure}
Just to give an idea, consider two concentric discrete spheres
of radii~$r<R$, in~$\Z^d$ or on a discrete torus.
Then, assume that one wants to study the \emph{trace}
left by a random walk (or random interlacements)
on the smaller ball. 
For this, consider the excursions between the two spheres
(we do not define the excursions formally, since
Figure~\ref{f_excursions} speaks for itself);
clearly, the trace of the whole process equals the trace of
these excursions.
Now, of course the excursions are not independent, since 
the law of the initial point of the next excursion 
depends on the last point of the previous excursion.
However, if we assume that $r\ll R$, this dependence
will be very weak, thus permitting one to use the results
of this paper for comparing these excursions to independent
excursions.

\section{Notations and results}
\label{notations}
We start describing the assumptions under which 
we will prove our main results.

Let~$(\Sigma,d)$ be a compact metric space, with~$\B(\Sigma)$ 
representing its Borel~$\sigma$-algebra.

\begin{assump}\label{assump_1}
 $(\Sigma,d)$ is of \emph{polynomial class}, that is, there 
exist some $\beta\geq 0$ and~$\phi\geq 1$ such that, 
for all $r\in(0,1]$, the number of open balls of radius at most~$r$
needed to cover $\Sigma$ is smaller than or equal
 to $\phi r^{-\beta}$. 
\end{assump}

As an example of metric space of polynomial class, consider first a
finite space~$\Sigma$, endowed with the discrete metric
\[
d(x,y)=\ind_{\{x\neq y\}}, \text{ for } x,y\in\Sigma.
\]
In this case, we can choose~$\beta=0$ and~$\phi=|\Sigma|$ (where~$|\Sigma|$ 
represents the cardinality of~$\Sigma$). As a second example, let us 
consider~$\Sigma$ to be a compact $k$-dimensional Lipschitz 
submanifold of~$\R^m$ with metric induced by the Euclidean norm of~$\R^m$. 
In this case we can take~$\beta=k$, but $\phi$ will in general depend 
on the precise structure of~$\Sigma$.
It is important to observe that, for a finite~$\Sigma$, it may not 
be the best idea to use the above discrete metric; one may be better off
with another one, e.g., the metric inherited from the Euclidean
space where~$\Sigma$ is immersed (see e.g.\ the proof 
of Lemma~2.9 of~\cite{CP}). 

We consider a Markov chain~$X=(X_i)_{i\geq 1}$ with transition
kernel~$\kP(x,dy)$ and starting law~$\V$, on~$(\Sigma,\B(\Sigma))$, and we suppose that
the chain has a unique invariant probability measure~$\Pi$. 
Moreover, we assume that the starting law and the transition kernel are absolutely continuous 
with respect to~$\Pi$. 
Let us denote respectively by~$\nu(\cdot)$ and~$p(x,\cdot)$ the Radon-Nikodym derivatives
(i.e., \emph{densities}) of~$\V(\cdot)$ and~$\kP(x,\cdot)$: for all~$A\in\B(\Sigma)$
\begin{linenomath}
\begin{align*}
\V(A) &= \int_{A} \nu(y) \Pi(dy), \\
\kP(x,A) &= \int_{A} p(x,y) \Pi(dy), \text{ for } x\in\Sigma.
\end{align*}
\end{linenomath}
We also use
\begin{assump}\label{assump_2}
The density $p(x,\cdot)$ is 
\emph{uniformly H\"older continuous}, that is,
there exist constants $\kappa>0$ 
and $\gamma\in(0,1]$ such that for all $x,z,z'\in \Sigma$,
\begin{linenomath}
\begin{equation*}
|p(x,z)-p(x,z')|\leq \kappa d^{\gamma}(z,z').
\end{equation*}
\end{linenomath}
\end{assump} 
We also work under
\begin{assump}\label{assump_3}
There exists $\eps\in (0,\frac{1}{2}]$ such that 
\begin{equation}
\label{max_eps}
\sup_{x,y\in\Sigma}|p(x,y)-1| \leq \varepsilon,
\end{equation}
and
\begin{equation}
\label{cond_nu}
\sup_{x\in \Sigma}|\nu(x)-1| \leq \eps.
\end{equation}
\end{assump}
\noindent
Observe that~\eqref{cond_nu} is not very restrictive because,
due to~\eqref{max_eps}, the chain will anyway
come quite close to stationarity already on step~$2$.

 Additionally, let us denote by $Y=(Y_i)_{i\geq 1}$ a sequence
of i.i.d.\ random variables with law~$\Pi$.

Before stating our main result, we recall the definition
 of the total 
variation distance between
two probability measures~$\bar{\mu}$ and~$\hat{\mu}$
 on some measurable 
space $(\Omega, \mathcal{T})$,
\begin{linenomath}
\begin{equation*}
\|\bar{\mu}-\hat{\mu}\|_{\text{TV}}
=\sup_{A\in\mathcal{T}}|\bar{\mu}(A)-\hat{\mu}(A)|.
\end{equation*}
\end{linenomath}
When dealing with random elements $U$ and $V$, we will write 
(with a slight abuse of notation) $\text{d}_{\text{TV}}(U,V)$ 
to denote the total variation distance between the laws of~$U$
 and~$V$.
Denoting by $L_n^Z:=(L_n^Z(x))_{x\in \Sigma}$ the local time 
field of the process $Z=X$ or $Y$ at time $n$, 
we are now ready to state 
{\thm \label{Main_Thm} Under Assumptions~\ref{assump_1}--\ref{assump_3}, 
there exists a universal positive constant $K$ such that, for all $n\geq 1$, it holds that
 \begin{linenomath}
\begin{align*}
\emph{d}_{\emph{TV}}(L_n^X , L_n^Y) \leq K\eps 
\displaystyle\sqrt{1 + \ln(\phi 2^{\beta}) 
+ \frac{\beta}{\gamma}\ln\Big(\frac{\kappa\vee (2\eps)}{\eps}\Big)}.
\end{align*} 
\end{linenomath}}

Note that the above bound is only useful when~$\eps$ is small enough;
however, we believe that, for such~$\eps$, this bound is
relatively sharp. 
As an application of Theorem~\ref{Main_Thm}, consider a 
finite state space~$\Sigma$, endowed with the discrete metric.
As we have already mentioned, in this case we can choose~$\beta=0$ 
and~$\phi=|\Sigma|$. Additionally, observe that, for any 
Markov chain~$X$ on $\Sigma$, under Assumption \ref{assump_3}, we can always take~$\kappa=2$ and~$\gamma=1$, 
so the uniform H\"older continuity of~$p$ 
is automatically verified here. Thus, Theorem~\ref{Main_Thm} leads to
\begin{linenomath} 
\begin{equation*}
\text{d}_{\text{TV}}(L_n^X , L_n^Y ) 
\leq K\eps \displaystyle\sqrt{1 + \ln|\Sigma|},
\end{equation*}
\end{linenomath}
for all~$n\geq 1$.

It is also relevant to check if we can obtain a uniform control 
(in $n$) of 
 $\text{d}_{\text{TV}}(L_n^X , L_n^Y)$ away from the 
``almost stationarity'' regime, i.e., 
for all~$\eps\in(0,1)$. In this direction, we introduce the following assumption, which is a slight modification of Assumption~\ref{assump_3}.
\begin{assump}\label{assump_4}
	There exists $\eps\in (0,1)$ such that 
	\begin{linenomath}
	\begin{equation*}
	\Big(\sup_{x,y\in\Sigma}|p(x,y)-1|\Big) \vee \Big(\sup_{x\in \Sigma}|\nu(x)-1|\Big) \leq \varepsilon.
	\end{equation*}
	\end{linenomath}
\end{assump}
We obtain the following
{\thm \label{Thm2} Under Assumptions~\ref{assump_1}, \ref{assump_2} and~\ref{assump_4}, there exists a positive 
cons\-tant~$K'=K'(\beta,\varphi,\kappa,\gamma,\eps)$, decreasing in $\eps$,
 such that
\begin{linenomath}
\begin{equation*}
\emph{d}_{\emph{TV}}(L_n^X , L_n^Y) \leq 1-K'
\end{equation*} 
\end{linenomath}
for all $n\geq 1$.
}

Such a result can be useful e.g.\ in the following context:
if we are able to prove that, for the i.i.d.\ sequence,
 something happens with probability close to~$1$, then 
the same happens for the field of local time of the Markov chain 
with at least uniformly positive probability. 
Observe that it is not unusual that the fact 
that the probability of something is uniformly positive
implies that it should be close to~$1$ then
(because one frequently has general results stating that 
this something should converge to~$0$ or~$1$).

Unlike Theorem~\ref{Main_Thm}, here we do not have 
a tractable explicit form for the above constant~$K'$; 
anyway, we think that with the method of the proof we use 
it is difficult to obtain a reasonably sharp expression for it
 (again, unlike Theorem~\ref{Main_Thm}).
Also, note that $\eps<1$ still means that the Markov chain
can regenerate in just one step with uniformly positive
probability. We conjecture that this can be relaxed
(by making a suitable assumption on e.g.\ the mixing time),
but, for now, we do not have a conclusive argument in this
direction.

The rest of the paper is organized in the following way.
In Section~\ref{Sim}, among other things, we show how the 
soft local time method  
 can be applied to the Markov chain~$X$ for constructing its 
 local time field. 
In Section~\ref{coupling} we present the construction of a 
coupling between the local time fields of 
the two processes~$X$ and~$Y$ at time~$n$. 
In Section~\ref{TV_binomial} we estimate the total variation 
distance between two binomial point processes. 
This auxiliary result
will be useful to bound from above the probability of the 
complement
of the coupling event introduced in Section~\ref{coupling}. 
In Section~\ref{Premres} we use a concentration inequality 
due to~\cite{Adam08} 
together with the machinery of empirical processes to obtain 
some intermediate results. In Section~\ref{Main_Thm_proof} 
we give the proof of Theorem~\ref{Main_Thm}. Finally, in Section~\ref{Second_Thm}, we give the proof of Theorem~\ref{Thm2}.

We end this section with considerations on the notation for 
constants used in this paper. Throughout the text, in general, we 
use capital~$C_1,C_2,\dots$ to denote global constants that 
appear in the results.
When these constants depend on some parameter(s), 
we will explicitly put (or mention) the dependence, otherwise the constants are 
considered universal. Moreover,  we will use small~$c_1,c_2,\dots$ 
to denote ``local constants'' that appear locally in the proofs, 
restarting the enumeration at the beginning of each proof.

\section{Constructions using soft local times}
\label{Sim}

We assume that the reader is familiar with the 
general idea of using Poisson point processes for constructing
general adapted stochastic processes, also known as
the \emph{soft local time} method. We refer to Section~4
of~\cite{PT15} for the general theory, and also
to Section~2 of~\cite{CGPV13} for a simplified introduction.

In this paper we use a modified version of
 this technique to couple the local time fields of
both processes, and we do this precisely in Section~\ref{coupling}. 
In this section 
 we first present two different constructions of the Markov chain~$X$, 
and then we present a construction of the local time field of~$X$ only. 
Then, in Section~\ref{Sim_IID}, we present a construction 
of the i.i.d.\ sequence~$Y$. All these constructions 
consist in applying the method of soft local times
in a (relatively) straightforward way.

Let~$\alpha$ be the \emph{regeneration coefficient} of the 
chain~$X$ with respect to~$\Pi$ (see Definition~4.28 of~\cite{FG}) defined by
\begin{linenomath}
\begin{align*}
\alpha := \inf_{x,y\in\Sigma} p(x,y).
\end{align*} 
\end{linenomath}
Note that $\alpha\leq 1$ since $p(x,\cdot)$ is a probability density 
for all $x\in \Sigma$.
Moreover, \eqref{max_eps} implies
that~$\alpha\geq 1-\eps\geq \frac{1}{2}$. 
 Hence, we consider the following decomposition 
 \begin{linenomath}
\begin{align*}
p(x,\cdot) = \frac{1}{2} + \frac{1}{2}\mu(x,\cdot), \text{ for all } x\in\Sigma,
\end{align*}
\end{linenomath}
where $\mu(x,\cdot)=2p(x,\cdot) -1\geq 0$ 
is a probability density with respect to~$\Pi$.

On some probability 
spa\-ce~$(\tilde{\Omega},\tilde{\mathcal{T}},\IP)$, suppose that we are 
given the following independent random elements:
\begin{itemize}
\item A sequence~$(I_j)_{j\geq 1}$ with~$I_1=1$ and~$(I_j)_{j\geq 2}$ 
i.i.d.~Bernoulli($\frac{1}{2}$) random variables;
\item A Poisson point process~$\eta$ on~$\Sigma\times\R_+$ 
with intensity measure~$\Pi\otimes\lambda_+$, where~$\lambda_+$
is the Lebesgue measure on~$\R_+$ and $\Pi$ is the invariant 
probability measure of~$X$ (cf.\ Section~\ref{notations}).  
\end{itemize}

Then, we define the sequence~$(\rho_j)_{j\geq 0}$ such that
\begin{linenomath}
\begin{align*}
\rho_0 &= 1,\\
\rho_{k+1} &= \inf\{ j>\rho_k : I_j=1 \} \text{ for } k\geq 0.
\end{align*}
\end{linenomath}
We interpret the elements of the sequence~$(\rho_j)_{j\geq 1}$ as being 
the random re\-ge\-ne\-ra\-tion times of the Markov 
chain~$X$: at each random 
time~$\rho_j$, the chain~$X$ starts afresh with law~$\Pi$. 
In this way, the chain will be viewed as a sequence of 
(independent) blocks (called {\it regeneration blocks}) 
with starting law~$\Pi$ 
and transitions 
according to~$\mu(\cdot,\cdot)$. 
Such blocks thus have lengths given by 
the differences of the subsequent elements of~$(\rho_j)_{j\geq 1}$.

The Poisson point process~$\eta$ will be used in the next 
sections to construct the local time fields of the Markov 
chain~$X$ and the i.i.d.\ sequence~$Y$ 
by means of the soft local time method.

\subsection{Construction of the local time field 
of the Markov chain~$X$}
\label{Sim_MC}
We first give two ways to construct the Markov chain~$X$ up to time~$n$
using soft local times. 
For that, 
we consider the Poisson point process described above,
\begin{linenomath}
\begin{align*}
\eta = \sum_{\lambda\in\Lambda} \bdelta_{(z_{\lambda}, t_{\lambda})},
\end{align*}
\end{linenomath}
(where $\Lambda$ is a countable index set),
and we proceed with the soft local time scheme in the classical way first. 
Denote by~$(x_i)_{i\geq 1}$ the elements of~$\Sigma$ 
which we will consecutively construct.

We begin with the construction of~$x_1$ by defining
\begin{linenomath}
\begin{align*}
\xi_1 &= \inf\big\{\ell\geq 0: \exists (z_{\lambda},t_{\lambda}) 
\text{ such that } \ell  \nu(z_{\lambda}) \geq t_{\lambda}\big\}, \\
G^X_1(x) &=  \xi_1   \nu(x), \;\text{for all}\; x\in \Sigma,
\end{align*}
\end{linenomath}
and~$(x_1,t_1)$ to be the unique pair~$(z_{\lambda},t_{\lambda})$ 
satisfying~$G^X_1(z_{\lambda}) = t_{\lambda}$. 

Then, once we have obtained the first state~$x_1$ visited by the chain~$X$, 
we proceed to the construction of the other ones. For~$i=2,3,\dots$, let
\begin{linenomath}
\begin{align*}
\xi_i &= \inf\big\{\ell\geq 0: \exists (z_{\lambda},t_{\lambda})
\notin\{(x_k,t_k)\}_{k=1}^{i-1} \text{ such that } \\
& \hspace{5cm} G^X_{i-1}(z_{\lambda}) + \ell  p(x_{i-1},z_{\lambda})
\geq t_{\lambda}\big\}, \\
G^X_i(x) &= G^X_{i-1}(x) + \xi_i   p(x_{i-1},x), \;\text{for all}\; x\in \Sigma,
\end{align*}
\end{linenomath}
and~$(x_i,t_i)$ to be the unique pair~$(z_{\lambda},t_{\lambda})$ out 
of the set~$\{(x_k,t_k)\}_{k=1}^{i-1}$ satisfying~$G^X_i(z_{\lambda}) 
= t_{\lambda}$.

Thus, after performing this iterative scheme for~$n$ iterations, 
we obtain the accumulated soft local time of the Markov chain~$X$ 
at time~$n$, which is given by
\begin{linenomath}
\begin{align*}
G^X_n(x) = \xi_1 \nu(x) + \sum_{k=2}^{n} \xi_k p(x_{k-1},x) ,
\end{align*}
\end{linenomath}
for~$x\in\Sigma$.

Next, we present an alternative construction of the same Markov chain 
taking into account the regeneration times of~$X$, using the Poisson point
process $\eta$ and the sequence~$(I_j)_{j\geq 1}$ of Bernoulli random variables introduced at 
the beginning of this section. Denote now by~$(\hat{x}_i)_{i\geq 1}$ 
the elements of~$\Sigma$ 
which we will consecutively construct in this alternative way. 
By construction, the sequence~$(\hat{x}_i)_{i\geq 1}$ will also
have the law of the Markov chain~$X$.

We begin with the construction of~$\hat{x}_1$ by defining
\begin{linenomath}
\begin{align*}
\hat{\xi}_1 &= \inf\big\{\ell\geq 0: \exists (z_{\lambda},t_{\lambda}) 
\text{ such that } \ell  \nu(z_{\lambda}) \geq t_{\lambda}\big\}, \\
\hat{G}^X_1(x) &=  \hat{\xi}_1   \nu(x), \;\text{for all}\; x\in \Sigma,
\end{align*}
\end{linenomath}
and~$(\hat{x}_1,\hat{t}_1)$ to be the unique pair~$(z_{\lambda},t_{\lambda})$ 
satisfying~$\hat{G}^X_1(z_{\lambda}) = t_{\lambda}$. 

Then, once we have obtained the first state~$\hat{x}_1$ 
visited by the chain~$X$, we proceed to the construction
of the other ones. For~$i=2,3,\dots$, define 
\begin{linenomath}
\begin{align*}
\hat{\xi}_i &= \inf\big\{\ell\geq 0: \exists (z_{\lambda},t_{\lambda})
\notin\{(\hat{x}_k,\hat{t}_k)\}_{k=1}^{i-1} \text{ such that } \\ 
& \hspace{4.5cm} \hat{G}^X_{i-1}(z_{\lambda}) 
+ \ell \big(I_i +(1-I_i)\mu(\hat{x}_{i-1},z_{\lambda}) \big)  \geq t_{\lambda}\big\}, \\
\hat{G}^X_i(x) &= \hat{G}^X_{i-1}(x) 
+ \hat{\xi}_i  \big(I_i +(1-I_i)\mu(\hat{x}_{i-1},x) \big),\;\text{for all}\; x\in \Sigma,
\end{align*}
\end{linenomath}
and~$(\hat{x}_i,\hat{t}_i)$ to be the unique 
pair~$(z_{\lambda},t_{\lambda})$ out of the 
set~$\{(\hat{x}_k,\hat{t}_k)\}_{k=1}^{i-1}$ 
satisfying~$\hat{G}^X_i(z_{\lambda}) = t_{\lambda}$.

Thus, after performing this iterative scheme for~$n$ iterations, 
we obtain the accumulated soft local time at 
time~$n$
\begin{linenomath}
\begin{align*}
\hat{G}^X_n(x) = \hat{\xi}_1 \nu(x) + \sum_{k=2}^{n} \hat{\xi}_k 
\big(I_k +(1-I_k)\mu(\hat{x}_{k-1},x) \big),
\end{align*}
\end{linenomath}
for~$x\in\Sigma$. 

Since in this paper we are interested in the random 
field of local times of the chain until time~$n$, 
the order of appearance of the states of~$X$ is not relevant 
for us and we will use the soft local time scheme in a slightly 
different way from that described above.
Specifically, we use the random variables $I_1,\dots, I_n$ as in 
the previous construction but now we first construct all the regeneration
blocks of size strictly greater than one and then the regeneration blocks
of size one.

We consider the Poisson point process~$\eta$ and 
the random variables~$I_1, I_2, \dots, I_n$. 
Then, we define the random set~$\mH\subset\{1,2,\dots,n\}$ as
\begin{linenomath}
\begin{align}
\mH= \big\{j\in\{2,3,\dots,n-1\} : I_j I_{j+1}=1 \} 
\cup \{j\in\{n\} : I_j=1\big\},
\label{set_H}
\end{align}
\end{linenomath} 
and the random permutation~$\mfS:\{1,2,\dots,n\}
\rightarrow\{1,2,\dots,n\}$ in the following way: 
\begin{itemize}
\item for~$j\in \mH^c$, $\mfS(j) = j-\sum_{i=2}^{j-1}
I_i I_{i+1}$,
\item for~$j\in \mH$, $\mfS(j) 
= |\mH^c| + |\{i\in\mH : i\leq j\}|$,
\end{itemize}
with the convention that~$\sum_{i=2}^{k} =0$ if~$k<2$.

Now, define 
\begin{linenomath}
\begin{align*}
\tilde{\xi}_1 &= \inf\big\{\ell\geq 0: \exists (z_{\lambda},t_{\lambda}) 
\text{ such that } \ell  \nu(z_{\lambda}) \geq t_{\lambda}\big\}, \\
\tilde{G}^X_1(x) &=  \tilde{\xi}_1   \nu(x),\;\text{for all}\; x\in \Sigma,
\end{align*}
\end{linenomath}
and~$(\tilde{x}_1,\tilde{t}_1)$ to be the unique 
pair~$(z_{\lambda},t_{\lambda})$ 
satisfying~$\tilde{G}^X_1(z_{\lambda}) = t_{\lambda}$. 
Next, for~$i=2,3,\dots,n$, define 
\begin{linenomath}
\begin{align*}
\tilde{\xi}_i &= \inf\big\{\ell\geq 0: \exists (z_{\lambda},t_{\lambda})
\notin\{(\tilde{x}_k,\tilde{t}_k)\}_{k=1}^{i-1} \text{ such that } \\ 
& \hspace{2cm} \tilde{G}^X_{i-1}(z_{\lambda}) 
+ \ell (I_{\mfS^{-1}(i)} + (1-I_{\mfS^{-1}(i)})
\mu(\tilde{x}_{i-1},z_{\lambda})) \geq t_{\lambda}\big\}, \\
\tilde{G}^X_i(x) &= \tilde{G}^X_{i-1}(x) + \tilde{\xi}_i 
(I_{\mfS^{-1}(i)} + (1-I_{\mfS^{-1}(i)}) 
\mu(\tilde{x}_{i-1},x)),\;\text{for all}\; x\in \Sigma,
\end{align*}
\end{linenomath}
and~$(\tilde{x}_i,\tilde{t}_i)$ to be the unique
pair~$(z_{\lambda},t_{\lambda})$ out of the 
set~$\{(\tilde{x}_k,\tilde{t}_k)\}_{k=1}^{i-1}$ 
satisfying~$\tilde{G}^X_i(z_{\lambda}) = t_{\lambda}$.

At the end of this procedure, we obtain the accumulated soft 
local time until time~$n$,
\begin{linenomath}
\begin{align*}
\tilde{G}^X_n(x) = \tilde{\xi}_1\nu(x) +\sum_{i=2}^{n} \tilde{\xi}_i  
(I_{\mfS^{-1}(i)} + (1-I_{\mfS^{-1}(i)}) 
\mu(\tilde{x}_{i-1},x)),
\end{align*}
\end{linenomath}
for~$x\in\Sigma$, observing that~$\tilde{G}^X_n$ has the same 
law as~$\hat{G}^X_n$, under~$\IP$. Also, observe 
that when proceeding in this way, we obtain the decomposition
\begin{linenomath}
\begin{align*}
\tilde{G}^X_n(x) = \tilde{G}^X_{|\mH^c|}(x) 
+ (\tilde{G}^X_n(x)-\tilde{G}^X_{|\mH^c|}(x)),
\end{align*}
\end{linenomath} 
where~$\tilde{G}^X_{|\mH^c|}(x)$ is the accumulated soft local time
corresponding to the construction of the first block plus the 
regeneration blocks of size strictly greater than one (until time~$n$) and $(\tilde{G}^X_n(x)-\tilde{G}^X_{|\mH^c|}(x))$ is the sum of the regeneration blocks of size one.

By implementing this last scheme, one produces a sequence~$(\tilde{x}_i)_{i}$ 
with~$n$ elements, the local time field of which has the law of the local time 
field of~$X$ at time~$n$, just as we wanted. Also, we recall the property
(of the soft local times) that the elements in the family~$(\tilde{\xi}_i)_i$ 
are all i.i.d.~Exponential($1$) random variables, independent 
of all other random elements (cf.~\cite{PT15}).

\subsection{Construction of the i.i.d.~sequence~$Y$}
\label{Sim_IID}
Now, we describe how we can construct the i.i.d.~sequence~$Y_1,\ldots,Y_n$ 
 using the soft local time technique. 

We denote by~$(y_i)_{i\geq 1}$ the elements of~$\Sigma$ which we 
will consecutively construct. 
Considering the same Poisson point process~$\eta$ described 
above, we begin with the construction of~$y_1$ by defining
\begin{linenomath}
\begin{align*}
\xi'_1 &= \inf\big\{\ell\geq 0: \exists (z_{\lambda},t_{\lambda}) 
\text{ such that } \ell \geq t_{\lambda}\big\}, \\
G^Y_1(x) &=  \xi'_1,\;\text{for all}\; x\in \Sigma,
\end{align*}
\end{linenomath}
and~$(y_1,t_1)$ to be the unique pair~$(z_{\lambda},t_{\lambda})$ 
satisfying~$G^Y_1(z_{\lambda}) = t_{\lambda}$.

Then, we proceed to the construction of~$y_2,y_3,\dots,y_n$: 
for~$i=2,3,\dots,n$, define 
\begin{linenomath}
\begin{align*}
\xi'_i &= \inf\big\{\ell\geq 0: \exists (z_{\lambda},t_{\lambda})
\notin\{(y_k,t_k)\}_{k=1}^{i-1} \text{ such that }G^Y_{i-1}(z_{\lambda}) 
+ \ell \geq t_{\lambda}\big\}, \\
G^Y_i(x) &= G^Y_{i-1}(x) + \xi'_i, \;\text{for all}\; x\in \Sigma,
\end{align*}
\end{linenomath}
and~$(y_i,t_i)$ to be the unique pair~$(z_{\lambda},t_{\lambda})$ 
out of the set~$\{(y_k,t_k)\}_{k=1}^{i-1}$ satisfying~$G^Y_i(z_{\lambda}) = t_{\lambda}$. 
At the end of this iterative scheme, 
we obtain the first~$n$ elements of the sequence~$Y$.   
As before, the elements in the family~$(\xi'_i)_i$ are all
i.i.d.\ Exponential($1$) random variables, and independent of 
all the other quantities.

The use of the soft local times to construct  the sequence~$Y$, 
as described above, produces the accumulated soft local time until time~$n$,
\begin{linenomath}
\begin{align*}
G_n^Y(x) =  \sum_{k=1}^{n} \xi'_k, \text{ for } x\in\Sigma.
\end{align*}
\end{linenomath}

\section{The coupling}
\label{coupling} 

In order to construct the coupling we are looking for, we assume that, 
in addition to the Bernoulli sequence 
$(I_j)_{j\geq 1}$ and the Poisson point 
process~$\eta$,  the triple $(\tilde{\Omega}, \tilde{\mathcal{T}}, \IP)$ 
(from Section~\ref{Sim}) also support all the other random elements to be introduced in this section. 

We start by partioning $\tilde{\Omega}$ using the event
\begin{linenomath}
\begin{align}
\label{event_C}
\mathsf{C}:=\Big\{|\mH|> \frac{n}{24}\Big\},
\end{align}
\end{linenomath}
(where~$\mH$ is defined in~\eqref{set_H}) and its complement $\mathsf{C}^c$. The idea is to use again the Bernoulli random variables $I_1,\dots, I_n$ to fix the regeneration times of the Markov chain $X$. In this way, the set $\mathsf{C}$ can be seen as the set of ``good realizations'' of $I_1,\dots, I_n$, in the sense that they produce a large number of regenerations for the chain~$X$. On the other hand, $\mathsf{C}^c$ can be seen as the set of ``bad realizations'' of $I_1,\dots, I_n$, in the sense that they produce few regenerations for $X$.

On $\mathsf{C}^c$, we will construct a rather rough coupling between the conditional law of $X$ (with few regenerations) and the i.i.d.~sequence $Y$, which will be enough for our purpose. On $\mathsf{C}$, we will construct a much more subtle coupling between the local times of the conditonal law of $X$ (with many regenerations) and $Y$ using the soft local time method.

We first present the construction of the coupling on $\mathsf{C}^c$. More precisely, on this set, we construct a point-by-point coupling between the Markov chain~$X$ and the i.i.d.\ sequence~$Y$, until time~$n$. For this, consider the random variables~$X'_1$ and~$Y'_1$ with respective laws~$\nu\text{d}\Pi$ and~$\Pi$, such that they are maximally coupled and the pair~$(X'_1,Y'_1)$ is independent of~$(I_1,\dots,I_n)$. Next, for~$k=2,3,\dots,n$, given~$(X'_{k-1},Y'_{k-1})$, we proceed as follows:
\begin{itemize}
	\item If~$I_k=1$, we consider the random variable~$X'_k$, independent of everything, with law~$\Pi$, and take~$Y'_k=X'_k$; 
	\item If~$I_k=0$, we consider~$(X'_k,Y'_k)$ such that~$X'_k$ is distributed according to~$\mu(X'_{k-1}, \cdot)\text{d}\Pi$, $Y'_k$ is distributed according to~$\Pi$, and~$X'_k$ and~$Y'_k$ are maximally coupled.
\end{itemize}

Now, we present the construction of a coupling between 
the local time fields of the Markov 
chain~$X$ and the i.i.d.\ sequence~$Y$ 
at time~$n$, on the set~$\mathsf{C}$. We mention that the new random elements to be introduced in the rest of this section will be independent of all the random elements already introduced. For the sake of brevity, we introduce the notation~$\I=(I_1,\dots,I_n)$. We will use the random element 
$\W:=(\I,\eta)$ and the auxiliary random elements~$V$, $V'$, $V''$ and~$\eta'$ 
(that we define later in this section),
 to construct a coupling between two point processes~$\eta_{X''}$ and~$\eta_{Y''}$ such that, under~$\IP$, both point processes are copies of~$\eta$. 
These copies will be such that
the third construction of Section~\ref{Sim_MC} 
applied to~$\eta_{X''}$ and the
construction of Section~\ref{Sim_IID} applied to~$\eta_{Y''}$ will 
 give high
probability of successful coupling of the local time fields of~$X$ and~$Y$, given~$\mathsf{C}$, for~$\eps$ sufficiently small.

\begin{figure}
\begin{center}
\includegraphics{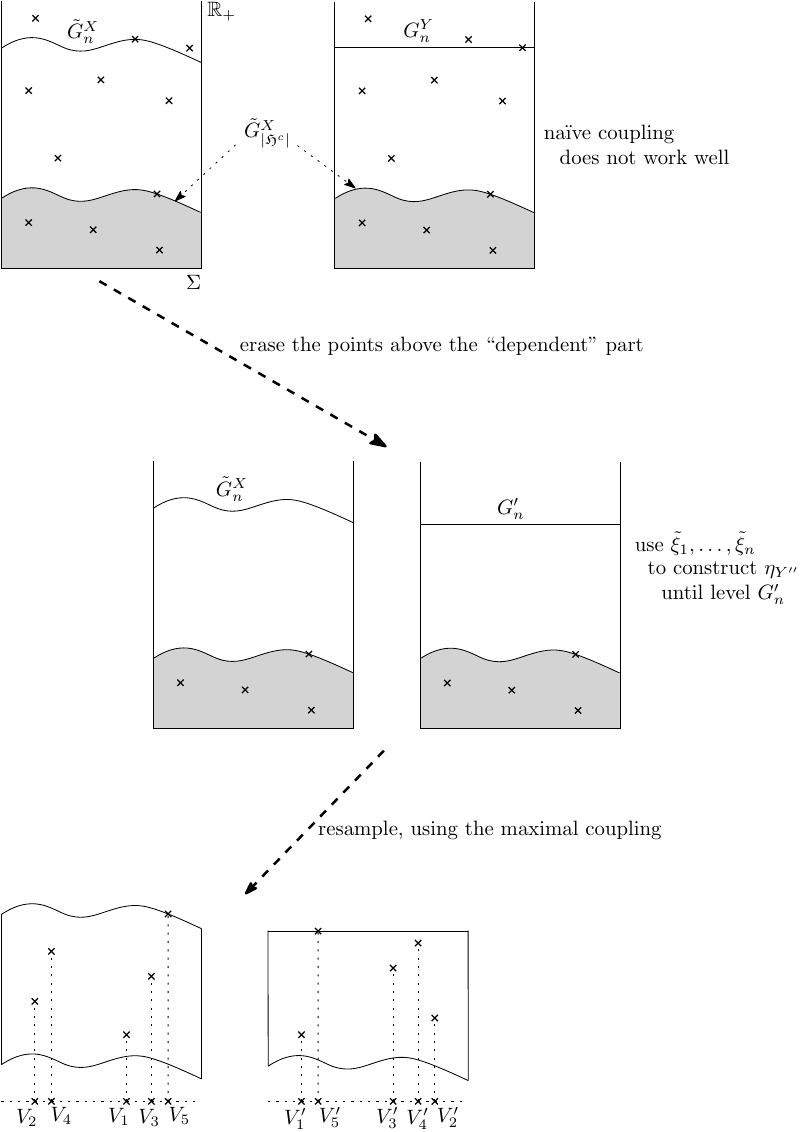}
\caption{Resampling of the ``independent parts'' on $\C$}
\label{f_resampling}
\end{center}
\end{figure}

It is important to stress that the ``na\"\i{}ve''
coupling (that is, using the same realization of the 
Poisson marks for constructing both the Markov chain
and the i.i.d.\ sequence) 
does not work, because it is not probable that both 
constructions will pick \emph{exactly} the same marks
 (look at the two marks at the upper 
right part of the top pictures on Figure~\ref{f_resampling}).
To circumvent this, we proceed as shown on Figure~\ref{f_resampling}:
 we first remove all the 
marks which are above the ``dependent'' part (that is, the marks strictly above the curve~$\tilde{G}^{X}_{|\mH^c|}$), and then resample
them using the maximal coupling of the ``projections''.
In the following, we describe this construction 
in a rigorous way.

We start with the construction of~$\eta_{X''}$ on $\mathsf{C}$.
As explained in Section~\ref{Sim_MC},
we first use~$\W$ to construct the local time field of 
the Markov chain~$X$ 
up to time~$n$. In this way, 
we obtain the soft local time curves
 $\tilde{G}^X_i$, for $1\leq i\leq n$, 
together with the sequences $\tilde{\xi}_1,\dots, \tilde{\xi}_n$
 and $\tilde{x}_1,\dots, \tilde{x}_n$.
 
Then, we define the random probability density (with respect to $\Pi$) $\Psi(\cdot)$ equal to
\begin{linenomath}
\begin{align}
 \frac{\Big(\displaystyle\sum_{j=|\mH^c|+1}^{n}\tilde{\xi}_j -\tilde{G}^{X}_{|\mH^c|}(\cdot)+\inf_{y\in \Sigma}\tilde{G}^{X}_{|\mH^c|}(y)\Big)_+}
{\displaystyle  \int_{\Sigma} \Big( \sum_{j=|\mH^c|+1}^{n}\tilde{\xi}_j- \tilde{G}^{X}_{|\mH^c|}(x)+\inf_{y\in \Sigma}\tilde{G}^{X}_{|\mH^c|}(y)\Big)_+\Pi(\text{d}x)},
\label{fct_psi}    
\end{align}
\end{linenomath}
if the denominator of the last expression is positive and $\mathds{1}(\cdot)$ otherwise.

 Anticipating on what is coming, on~$\mathsf{C}$, the law
 $\Psi \text{d}\Pi$ will be used to sample 
the first coordinates of the~$|\mH|$ marks of~$\eta_{Y''}$ between $\tilde{G}^{X}_{|\mH^c|+1}$
 and~$G'_n:=\inf_{y\in \Sigma}\tilde{G}^{X}_{|\mH^c|}(y)+\sum_{j=|\mH^c|+1}^{n}\tilde{\xi}_j$. 

Now, going back to the construction of~$\eta_{X''}$, we adopt a resampling scheme: 
we first ``erase'' 
all the marks of the point process~$\eta$ in the space~$\Sigma\times\R_+$ that are on the curves 
$\tilde{G}^X_{|\mH^c|+1},\dots,\tilde{G}^X_n$, 
then we reconstruct the marks as follows.
We introduce the random vector $V:=(V_1,\dots,V_{_{|\mH|}})$
 such that under $\IP[\;\cdot \mid \W]$, its coordinates are 
independent and distributed according to the invariant measure~$\Pi$.
 We use the random vector~$V$ to place the (new) marks
\[
 \Big(V_1,\tilde{G}^X_{|\mH^c|+1}(V_1)\Big),
\Big(V_2,\tilde{G}^X_{|\mH^c|+2}(V_2)\Big),\dots, 
\Big(V_{|\mH|},\tilde{G}^X_{n}(V_{|\mH|})\Big),
\]
  on the curves $\tilde{G}^X_{|\mH^c|+1},\dots,\tilde{G}^X_n$ (see Figure~\ref{f_resampling}, bottom left picture).
 Finally, to obtain a point process $\eta_{X''}$ well defined on all $\tilde{\Omega}$, we complete its construction by defining $\eta_{X''}:=\eta$ on $\C^c$.

We continue with the construction of $\eta_{Y''}$. As we have just done for $\eta_{X''}$, we can get rid of the construction of $\eta_{Y''}$ on $\C^c$, by stating $\eta_{Y''}:=\eta$.

On $\C$, we will first construct the marks of $\eta_{Y''}$ between $\tilde{G}^X_{|\mH^c|}$ and $G'_n$. 
First, consider the random vector 
$V':=(V'_1,\dots,V'_{|\mH|})$ and the random element $V''=(V''_1(x),\dots ,V''_n(x))_{x\in \Sigma}$ such that 
under $\IP[\;\cdot \mid \W=w]$, $V'$ and $V''$ are independent, and for $w\in \tilde{\Omega}$,
\begin{itemize}
\item $V'_1,\dots,V'_{|\mH|}$ are independent;
\item $V'_{|\mH|}$ has law $\frac{\1{G'_n>\tilde{G}^X_{|\mH^c|}}}{\Pi[G'_n>\tilde{G}^X_{|\mH^c|}]} \text{d}\Pi$ if $\Pi[G'_n>\tilde{G}^X_{|\mH^c|}]>0$, and law $\Pi$ otherwise. Also, $V'_{|\mH|}$ is maximally coupled with $V_{|\mH|}$;
\item the elements $V'_1,\dots,V'_{|\mH|-1}$ have law $\Psi \text{d}\Pi$ and the elements\\
 $\big(\sum_{i=1}^{|\mH|-1}
\ind_{\{x\}}(V'_i)\big)_{x\in  \Sigma}$ 
and 
$\big(\sum_{i=1}^{|\mH|-1}\ind_{\{x\}}(V_i)\big)_{x\in  \Sigma }$ 
are maximally coupled;
\item the elements $(V''_1(x))_{x\in \Sigma},\dots,(V''_n(x))_{x\in \Sigma}$ are i.i.d.;
\item $(V''_1(x))_{x\in \Sigma}$ is a family of independent random variables and, for $x\in \Sigma$, $V''_1(x)$ has law $U(0, G'_n(x)-\tilde{G}^X_{|\mH^c|}(x))$ if $\Psi(x)>0$ and law $U(0,1)$ if $\Psi(x)=0$.
\end{itemize}

We construct the marks of the point 
process~$\eta_{Y''}$ below $G'_n$ in the following way:
we keep the marks obtained below $\tilde{G}^{X}_{|\mH^c|}$ 
and we use the law $\Psi\text{d}\Pi$ to complete the process 
until~$G'_n$.
For this, we adopt a resampling scheme just as before. 
We first erase all the marks of the point process~$\eta$
 that are (strictly) above $\tilde{G}^{X}_{|\mH^c|}$, 
then we resample the part of the process $\eta$ up to~$G'_n$, 
using the marks:
\begin{linenomath}
\begin{align*}
\lefteqn{\Big(V'_1,\tilde{G}^X_{|\mH^c|}(V'_1)+V''_1(V'_1)\Big),\dots, \Big(V'_j,\tilde{G}^X_{|\mH^c|}(V'_j)+V''_j(V'_j)\Big),\dots,}\nonumber\\
&\phantom{******}\Big(V'_{|\mH|-1},\tilde{G}^X_{|\mH^c|}(V'_{|\mH|-1})+V''_{|\mH|-1}(V'_{|\mH|-1})\Big),\Big(V'_{|\mH|},G'_n(V'_{|\mH|})\Big)
\end{align*}
\end{linenomath}
(see Figure~\ref{f_resampling}, bottom right picture). 

 
 Finally, we consider a copy $\eta'=\sum_{\lambda\in \Lambda}\bdelta_{(z'_{\lambda}, t'_{\lambda})}$ of~$\eta$, independent of the all the other random elements already introduced. Then, we use the marks of the point process $T\eta':=\sum_{\lambda\in \Lambda}\bdelta_{(z'_{\lambda}, t'_{\lambda}+G'_n\vee \tilde{G}^X_n(z'_{\lambda}))}$ to complete the marks of~$\eta_{Y''}$ above~$G'_n\vee \tilde{G}^X_n $, on~$\C$.
 
For the sake of brevity, let us denote by $m^{X''}_1,\dots, m^{X''}_n$ and $m^{Y''}_1,\dots, m^{Y''}_n$, the first coordinates ($\in \Sigma$) of the marks of $\eta_{X''}$ and $\eta_{Y''}$ below the curves
$\tilde{G}^X_n$ and~$G^Y_n$ (see Section \ref{Sim_IID}) respectively. Let~$L^{X''}_n$ and~$L^{Y''}_n$ be the fields of
 local times associated to these first coordinates, that is, for all $x\in \Sigma$,
 \begin{linenomath}
\begin{equation*}
L^{X''}_n(x)=\sum_{i=1}^n\ind_{\{x\}}(m^{X''}_i)\phantom{**} \text{and}\phantom{**} L^{Y''}_n(x)=\sum_{i=1}^n\ind_{\{x\}}(m^{Y''}_i).
\end{equation*}
\end{linenomath}
By construction, we have the following
 {\prop 
It holds that
\begin{linenomath}
\begin{align*}
L^{X}_n &\eqlaw \ind_{\C^c}L^{X'}_n + \ind_{\C}L^{X''}_n, \\
L^{Y}_n &\eqlaw \ind_{\C^c}L^{Y'}_n + \ind_{\C}L^{Y''}_n,
\end{align*}
\end{linenomath}
where $\eqlaw $ stands for equality in law. 
}
\medskip

\noindent
Consequently, we obtain a coupling between~$L^X_n$ and~$L^Y_n$. 
We will denote by~$\Upsilon$ the coupling event associated to 
this coupling, that is, 
\begin{linenomath}
\begin{align*}
\Upsilon = \Big\{\ind_{\C^c}L^{X'}_n + \ind_{\C}L^{X''}_n=\ind_{\C^c}L^{Y'}_n + \ind_{\C}L^{Y''}_n\Big\}.
\end{align*}
\end{linenomath}
 In Section~\ref{Main_Thm_proof}, 
we will obtain an upper bound for $\IP[\Upsilon^c]$.

\section{Total variation distance between binomial point processes} \label{TV_binomial}

In this section, we estimate the total variation distance 
between two binomial point processes on some measurable 
space $(\Omega, \mathcal{T})$ with laws~${\P}_n$ and~${\Q}_n$ 
of respective parameters $({\bf p}_n,n)$ and $(\q_n,n)$, 
where $n\in \mathbb{N}$ and $\p_n$, $\q_n$ are two probability 
laws on $(\Omega, \mathcal{T})$. We also assume 
that $\q_n\ll \p_n$ and that~$\p_n$ and~$\q_n$ 
are close in a certain sense to be defined below. 

For two probability measures~$\bar{\mu}$ and~$\hat{\mu}$ 
on $(\Omega, \mathcal{T})$, we recall that 
if $\bar{\mu}\ll \hat{\mu}$,
\begin{equation}
\label{TVdiscrete}
\|\bar{\mu}-\hat{\mu}\|_{\text{TV}}
=\frac{1}{2}\int_{\Omega}
\Big|\frac{\text{d}\bar{\mu}}{\text{d}\hat{\mu}}-1\Big|d\hat{\mu}.
\end{equation}

We will prove the following result, which is actually 
a little bit more than we need in this paper. 
{\prop 
\label{Propmulti}
Let $\delta_0\in (0,1]$ and $\delta \in [0,\delta_0)$ 
such that for all $n\in \mathbb{N}$, 
$|\frac{\emph{d}\q_n}{\emph{d}\p_n}(x)-1|\leq\delta n^{-1/2}$ 
for all~$x\in \Omega$. Then, for
 $C_1(\delta_0)=\exp(\delta_0^2)
\frac{\sinh(\delta_0^2)}{\delta_0}+\sqrt{2\pi}\exp(\frac{5}{2}\delta_0^2)$ 
we have, for all $n\in \mathbb{N}$,
\begin{linenomath}
\begin{equation*}
\|\P_n-\Q_n\|_{\emph{TV}}\leq C_1(\delta_0)\delta.
\end{equation*}
\end{linenomath}

}

\begin{proof}
In this proof, when we want to emphasize the probability 
law under which we take the expectation we will indicate the 
law as a subscript. For example, the expectation under some 
probability law~$\bar{\mu}$ will be denoted by~$E_{\bar{\mu}}$.

To begin, let us suppose that $n\geq 2$. We first observe 
that~${\P}_n$ and~${\Q}_n$ can be seen as probability measures 
on the space of $n$-point measures $\M_n=\{m:m=\sum_{i=1}^{n}\bdelta_{x_i}, 
x_i\in \Omega, 1\leq i\leq n\}$ endowed with 
the $\sigma$-algebra generated by the 
mappings $\Phi_B:\M_n\to \Z_+$ defined by $\Phi_B(m)=m(B)=\sum_{i=1}^n\bdelta_{x_i}(B)$, 
for all $B\in \mathcal{T}$. Observe that the law of~$\P_n$ 
(respectively,~$\Q_n$) is completely characterized by its values 
on the sets of the form $\{m\in\M_n:m(B_1)=n_1,\dots, m(B_J)=n_J\}$,
 where $J\in \Z_+$, $B_1,\dots, B_J$ are disjoint sets
 in~$\mathcal{T}$ and $n_1,\dots,n_J$ are non-negative integers
 such that $n_1+\dots+n_J=n$.
With this observation it is easy to deduce that $\Q_n\ll \P_n$ 
and check that its Radon-Nikodym derivative with respect 
to~$\P_n$ is given by
\[
\frac{\text{d}\Q_n}{\text{d}\P_n}(m)
=\prod_{i=1}^n\frac{\text{d}\q_n}{\text{d}\p_n}(x_i)
\]
where $m=\sum_{i=1}^{n}\bdelta_{x_i}$.\\
By (\ref{TVdiscrete}) we obtain that
\begin{linenomath}
\begin{align*}
\|\P_n-\Q_n\|_{\text{TV}}&=\frac{1}{2}\int_{\M_n}\Big|
\frac{\text{d}\Q_n}{\text{d}\P_n}(m)-1\Big|\text{d}\P_n(m)\nonumber\\
&=\frac{1}{2}\int_{\M_n}\Big|\prod_{i=1}^n
\frac{\text{d}\q_n}{\text{d}\p_n}(x_i)-1\Big|\text{d}\P_n(m).
\end{align*}
\end{linenomath}

Now, for all $n\in \mathbb{N}$, we define the function 
$f_n:\Omega\to \R$ such that, for~$x\in \Omega$, 
we have $f_n(x)=\frac{\text{d}\q_n}{\text{d}\p_n}(x)-1$. 
Observe that $E_{\p_n}[f_n]=0$ and that 
$\|f_n\|_{\infty}\leq\delta n^{-1/2}$ for all $n\geq 2$. We have that
\begin{linenomath}
\begin{align}
\label{TV1}
\|\P_n-\Q_n\|_{\text{TV}}&=\frac{1}{2}\int_{\M_n}\Big|
\prod_{i=1}^n(1+f_n(x_i))-1\Big|\text{d}\P_n(m)\nonumber\\
&=\frac{1}{2}\int_{\M_n}\Big|
\exp\big(\sum_{i=1}^ng_n(x_i)\big)-1\Big|
\text{d}\P_n(m)\nonumber\\
&=\frac{1}{2}E|\exp\{g_n(X_1)+\dots+g_n(X_n)\}-1|
\end{align}
\end{linenomath}
where, for all $n\geq 2$, $g_n$ is the function $\Omega\to \R$
 defined by $g_n=\ln(1+f_n)$ 
(by using $|\frac{\text{d}\q_n}{\text{d}\p_n}(x)-1|\leq 
\delta n^{-1/2}$ for all~$x\in \Omega$,
we can observe that $g_n$ is well defined) 
and, in the last equality, the random variables $X_1,\dots,X_n$ are i.i.d.\ with law~$\p_n$. 

Using the fact that $|\ln(1+x)|\leq 2|x|$ 
for $x\in (-1/\sqrt{2},1/\sqrt{2})$, we deduce that 
$\|g_n\|_{\infty}\leq 2\|f_n\|_{\infty}\leq 2\delta n^{-1/2}$, 
for all $n\geq 2$.
Observe that $|E[g_n(X_1)]|=|E[(g_n-f_n)(X_1)]|$ 
since $E[f_n(X_1)]=E_{\p_n}[f_n]=0$. 
Now we use the fact that, for all $x\in \R$ such 
that $|x|\leq 1/\sqrt{2}$, we have that
\[
|\ln(1+x)-x|\leq 2x^2.
\]
Since $\|f_n\|_{\infty} \leq \delta n^{-1/2}$ we obtain that 
$\|g_n-f_n\|_{\infty}\leq 2\delta^2 n^{-1}$. We deduce that 
\begin{linenomath}
\begin{align}
\label{EST1}
|E[g_n(X_1)]|\leq \|g_n-f_n\|_{\infty}\leq 2\delta^2 n^{-1}.
\end{align}
\end{linenomath}

Let $Z_n:=\sum_{k=1}^{n}f_n(X_k)$.  Using the fact that 
$|\exp(x) -1|\leq \exp(|x|)-1$ for all $x\in \R$ and~(\ref{EST1})
 we obtain that
 \begin{linenomath}
\begin{align}
\label{TV2}
\|\P_n-\Q_n\|_{\text{TV}} &\leq \frac{1}{2}E\Big[\exp\Big(|Z_n + \sum_{k=1}^{n}(g_n-f_n)(X_k)|\Big)-1\Big]\nonumber\\
&\leq \frac{1}{2}\Big(\exp(2\delta^2) E\Big[\exp(|Z_n|)\Big]-1\Big).
\end{align}
\end{linenomath}

Now let us obtain an upper bound for the expectation of the 
right-hand side of~(\ref{TV2}). Using integration by 
parts, we have
\begin{linenomath}
\begin{align}
\label{Parts}
E\Big[\exp(|Z_n|)\Big]=1+\int_{0}^{\infty}e^tP[|Z_n|
\geq t]dt.
\end{align}
\end{linenomath}
By Hoeffding's inequality, we have, for all $n\geq 2$ and 
for all $t\geq 0$,
\begin{linenomath}
\begin{equation*}
P[|Z_n|\geq t]\leq 2\exp\Big(-\frac{t^2}{2\delta^2}\Big).
\end{equation*}
\end{linenomath} 
Using this last inequality in~(\ref{Parts}), we obtain that
\begin{linenomath}
\begin{align*}
E\Big[\exp(|Z_n|)\Big]&=1+2\int_{0}^{\infty}
\exp\Big(t-\frac{t^2}{2\delta^2}\Big)dt\nonumber\\
&\leq 1+2\delta \exp\Big(\frac{\delta^2}{2}\Big)\sqrt{2\pi}.
\end{align*}
\end{linenomath}

Therefore, going back to~(\ref{TV2}) and using the fact that 
$\delta< \delta_0$, we deduce that for all $n\geq 2$,
\begin{linenomath}
\begin{align}
\label{EST3}
\|\P_n-\Q_n\|_{\text{TV}} &\leq \frac{1}{2}\Big(\exp(2\delta^2)
\Big[1+2\delta \exp\Big(\frac{\delta^2}{2}\Big)\sqrt{2\pi}\Big]-1\Big)\nonumber\\
&= \frac{1}{2} \exp(\delta^2)\Big[\exp(\delta^2)-\exp(-\delta^2)+2\delta\sqrt{2\pi}\exp\Big(\frac{3\delta^2}{2}\Big)\Big] \nonumber\\
&= \delta\Big[\exp(\delta^2)
\frac{\sinh(\delta^2)}{\delta}+\sqrt{2\pi}\exp\Big(\frac{5}{2}\delta^2\Big)\Big]\nonumber\\
&\leq c_1(\delta_0) \delta,
\end{align}
\end{linenomath} 
where $c_1(\delta_0):=\exp(\delta_0^2)
\frac{\sinh(\delta_0^2)}{\delta_0}+\sqrt{2\pi}\exp(\frac{5}{2}\delta_0^2)$.

Gathering~(\ref{TV1}), (\ref{EST3}) and considering the fact that
 $\|\P_1-\Q_1\|_{\text{TV}}\leq \frac{\delta}{2}$, 
we conclude the proof of Proposition~\ref{Propmulti}.
\end{proof}

\section{Controlling the ``dependent part'' of the soft
local time} 
\label{Premres}
First, recalling the notations of Section \ref{Sim}, we define the random function
\begin{linenomath}
	\begin{equation}
	\label{widehatW}
	\widehat{W}(x):= \sum_{k=\rho_1}^{\rho_2-1} \hat{\xi}_k(1-I_k)(1-\mu(\hat{x}_{k-1},x)),
	\;\text{for all}\; x\in \Sigma,
	\end{equation}
\end{linenomath}
Then, we consider a sequence of i.i.d.\ random functions 
$(\widehat{W}_n)_{n\geq 1}$ with the same law as $\widehat{W}$. 
We will show in Section~\ref{Est_F} that 
\begin{linenomath}
\begin{equation}
\label{EPE}
\sup_{n\in\IN}\frac{1}{\sqrt{n}}E\Big[\sup_{x\in \Sigma}\Big|\sum_{k=1}^{n}\widehat{W}_k(x)\Big|
\Big]\leq F
\end{equation}
\end{linenomath}
where
\begin{linenomath}
\begin{align}
	F := C_4\eps \displaystyle\sqrt{1 + \ln(\phi 2^{\beta}) + \frac{\beta}{\gamma}\ln\Big(\frac{\kappa\vee (2\eps)}{\eps}\Big)}
	\label{F_exp}
\end{align}
\end{linenomath}
and $C_4$ is a universal positive constant.

For all~$i\in\IN$ we define the events
\begin{linenomath}
\begin{equation}
\A_i=\Big\{1\wedge \sup_{x\in \Sigma}|\Psi(x)-1|\leq \frac{(1+i)F}{\sqrt{n}} \Big\}.
\label{event_Av}
\end{equation}
\end{linenomath}



The goal of this section is to prove the following
{\prop  \label{goodenv} There exist a universal positive constant~$C_2$ such that, for all~$n\geq 96$ and
 $i\in\IN$, it holds that
 \begin{linenomath}
\begin{align*}
\IP[\A_i^c\; |\; \C ] &\leq  \frac{C_2}{(1+i)^3}.
\end{align*}
\end{linenomath} 
}

We postpone the proof of this proposition to 
Section~\ref{Proof_Av}. Before that, in Section~\ref{Est_F} we show~\eqref{EPE} and in Section~\ref{concentr} 
we use a concentration inequality to obtain a tail estimate on some random variable related to $\Psi-1$ (see~Proposition \ref{Prop_Rn_tilde} below).

\subsection{Proof of inequality~\eqref{EPE}}
\label{Est_F}

In this section, we present a standard method based on bracketing numbers to prove~(\ref{EPE}). 
Without loss of generality, we assume 
in this section that~$\kappa$ from 
Assumption~\ref{assump_2} is greater than or equal to~$2\eps$.

We start introducing 
the space~$\mS=\R^{\Sigma}$ and the class~$\F=(f_x)_{x\in\Sigma}$ of 
functions~$f_x:\mS\rightarrow \R$ such that~$f_x(\omega)=\omega(x)$, 
for~$\omega\in\mS$. A function~$\mE:\mS\rightarrow\R$ is an \textit{envelope function} of the class~$\F$ if~$\mE\geq |f_x|$ for all~$f_x\in\F$.

In this setting, for~$s>0$, we need to estimate the bracketing number
 \begin{linenomath}
\begin{align*}
N_{[\,]}\Big(s\|\mE\|_{2},\F,L_2\Big), 
\end{align*}
\end{linenomath}
which is defined to be the minimum number of brackets
\begin{linenomath}	
\begin{align*}
[f_1,f_2] := \big\{f:\mS\rightarrow \R; f_1\leq f\leq f_2\big\},
\end{align*}
\end{linenomath}
satisfying~$\|f_2-f_1\|_{2} < s\|\mE\|_{2}$, that are
 needed to cover the class~$\F$, where the given 
functions~$f_1$ and~$f_2$ have finite $L_2$-norms 
(see Definition~2.1.6 of~\cite{VW96}).

For that, we consider an (initially arbitrary) exhaustive and 
finite collection of subsets of the space~$\Sigma$, that is, 
a finite collection~$\{\D_i\}_i$ such that~$\D_i\subset\Sigma$ 
for each~$i$ and~$\bigcup_i\D_i=\Sigma$, and for each such set~$\D_i$ 
we define two functions~$f_{\D_i},\hat{f}_{\D_i}:\mS\rightarrow\R$,
\begin{linenomath}
\begin{align*}
f_{\D_i}(\omega) = \inf_{z_0\in\D_i} \omega(z_0) ~\mbox{ and }~ \hat{f}_{\D_i}(\omega) = \sup_{z_0\in\D_i} \omega(z_0),
\end{align*}
\end{linenomath}
so that, if~$x\in\D_i$ then~$f_x\in[f_{\D_i},\hat{f}_{\D_i}]$.
 Thus, to each set in the family~$\{\D_i\}_i$ we associate 
a bracket, and so any particular finite exhaustive family of subsets
 of~$\Sigma$ induces a finite collection of 
brackets~$\{[f_{\D_i},\hat{f}_{\D_i}]\}_i$ which cover the class~$\F$.

So, in order to properly estimate the bracketing number, 
the task is to determine a suitable collection~$\{\D_i\}_i$ 
of subsets of~$\Sigma$ in such a way that the induced brackets have their
 sizes all smaller than~$s\|\mE\|_{2}$. The number of sets in that
 collection will serve as an upper bound for~$N_{[\,]}$.

Through the following lemma we better characterize the size
 (in~$L_2$) of the induced brackets we are considering.

{\lem Under Assumption~\ref{assump_2}, it holds that, for any
 set~$\D\subset\Sigma$,
 \begin{linenomath}
\begin{align*}
\|\hat{f}_{\D}-f_{\D}\|_{2} \leq 4\sqrt{2}\kappa
\displaystyle\max_{z,z'\in\D}d^{\gamma}(z,z').
\end{align*}
\end{linenomath}
\label{bracket_size}
}

\begin{proof}
Recalling the notation introduced at the beginning of Section~\ref{coupling}, 
we want to bound the~$L_2$-norm 
\begin{linenomath}
\begin{align*}
\Big\|\hat{f}_{\D}\Big(\widehat{W}_1(\cdot)\Big)-f_{\D}\Big(\widehat{W}_1(\cdot)
\Big)\Big\|_{2} 
= \Big\|\sup_{z\in\D}\widehat{W}_1(z)-\inf_{z\in\D}\widehat{W}_1(z)\Big\|_{2}
\end{align*}
\end{linenomath}
which is
\begin{linenomath}
\begin{align*}
\Big\| \sup_{z\in\D} \sum_{k=\rho_1}^{\rho_2-1} &\hat{\xi}_k(1-I_k)(1-\mu(\hat{x}_{k-1},z)) 
- \inf_{z\in\D} \sum_{k=\rho_1}^{\rho_2-1} \hat{\xi}_k(1-I_k)(1-\mu(\hat{x}_{k-1},z))\Big\|_{2}\\
&= \Big\| \sup_{z\in\D} \sum_{k=\rho_1}^{\rho_2-1} \hat{\xi}_k (1-I_k) \mu(\hat{x}_{k-1},z) -
 \inf_{z\in\D} \sum_{k=\rho_1}^{\rho_2-1} \hat{\xi}_k (1-I_k) \mu(\hat{x}_{k-1},z)\Big\|_{2} \\
&\leq \Big\| \sum_{k=\rho_1}^{\rho_2-1} \hat{\xi}_k (1-I_k) \Big(\sup_{z\in\D}  \mu(\hat{x}_{k-1},z) 
-  \inf_{z\in\D}  \mu(\hat{x}_{k-1},z)\Big) \Big\|_{2} \\
&\leq 2\kappa\Big\| \sum_{k=\rho_1}^{\rho_2-1} \hat{\xi}_k \Big\|_{2}  \displaystyle\max_{z,z'\in\D}d^{\gamma}(z,z'),
\end{align*}
\end{linenomath}
where we used Assumption~\ref{assump_2} and the definition of~$\mu$ to establish
 the second inequality.

Thus, we conclude the proof by using the fact 
that~$\sum_{k=\rho_1}^{\rho_2-1} \hat{\xi}_k$ is exponentially distributed
 with parameter~$\frac{1}{2}$, so that
 \begin{linenomath}
\begin{align*}
\Big\| \sum_{k=\rho_1}^{\rho_2-1} \hat{\xi}_k \Big\|_{2} = 2\sqrt{2}.
\end{align*}
\end{linenomath}
\end{proof}

In view of the above result, 
we should impose the sets we are constructing, 
$\{\D_i\}_i$, to be such that
\begin{linenomath}
\begin{align}
\max_{z,z'\in\D_i}d(z,z') < \Big(\frac{s\|\mE\|_{2}}{4\sqrt{2}\kappa}\Big)^{1/\gamma}, \mbox{ for each } i,
\label{condition_max}
\end{align}
\end{linenomath}
in order to obtain, from Lemma~\ref{bracket_size}, that
\begin{linenomath}
\begin{align}
\|\hat{f}_{\D_i}-f_{\D_i}\|_{2} < s\|\mE\|_{2}, \mbox{ for each } i.
\label{condition_max_2}
\end{align}
\end{linenomath}

Then, we prove 

{\prop Under Assumptions~\ref{assump_1} and~\ref{assump_2}, 
there exists a universal positive constant~$C_3$ such that, if
\begin{linenomath}
\begin{align}
\frac{\gamma}{\beta} \Big[1 + \ln(\phi 2^{\beta}) 
+ \frac{\beta}{\gamma}\ln\Big(\frac{4\sqrt{2}\kappa}
{\|\mE\|_{2}}\Big)\Big]  \geq \frac{1}{2},
\label{cond_geq1}
\end{align}
\end{linenomath}
then it holds that, for all~$n\in\IN$,
\begin{linenomath}
\begin{align*}
E\Big[\sup_{x\in \Sigma}\Big|\sum_{k=1}^{n}\widehat{W}_k(x)\Big|
\Big] \leq C_3 \displaystyle\sqrt{1 
+ \ln(\phi 2^{\beta}) + \frac{\beta}{\gamma}
\ln\Big(\frac{4\sqrt{2}\kappa}{\|\mE\|_{2}}\Big)} \|\mE\|_{2}\sqrt{n}.
\end{align*}
\end{linenomath}
\label{expect_R_tilde}
}

\begin{proof}
For~$s\in(0,1)$, we just define the sets~$\{\D_i\}_i$ to
 be open balls of radius at most
 \begin{linenomath}
\begin{align*}
\frac{1}{2}\Big(\frac{s\|\mE\|_{2}}{4\sqrt{2}\kappa}\Big)^{1/\gamma}<1,
\end{align*}
\end{linenomath}
so that~\eqref{condition_max} is verified (and
 consequently~\eqref{condition_max_2} too, 
by Lemma~\ref{bracket_size}).
Then, under Assumption~\ref{assump_1}, we can say that 
the number of such balls that are needed to cover~$\Sigma$ 
is at most
\begin{linenomath}
\begin{align*}
\phi 2^{\beta}\Big[\frac{4\sqrt{2}\kappa}
{s\|\mE\|_{2}}\Big]^{\beta/\gamma},
\end{align*}
\end{linenomath}
and therefore, for any~$s\in(0,1)$,
\begin{linenomath}
\begin{align*}
N_{[\,]}\Big(s\|\mE\|_{2},\F,L_2\Big) 
\leq \phi 2^{\beta}\Big[\frac{4\sqrt{2}\kappa}
{s\|\mE\|_{2}}\Big]^{\beta/\gamma}.
\end{align*}
\end{linenomath}

Taking this bound on the bracketing number into account,
 we can estimate the bracketing entropy integral (of the class~$\F$)
 \begin{linenomath}
\begin{align*}
J_{[\,]}\Big(1,\F,L_2\Big) := \int_{0}^{1} 
\sqrt{1+\ln N_{[\,]}\Big(s\|\mE\|_{2},\F,L_2\Big)} ds,
\end{align*}
\end{linenomath}
(see its definition e.g.\ in Section~2.14.1 
of~\cite{VW96}, page 240). We do that by just bounding it from above by
\begin{linenomath} 
\begin{align*}
\int_{0}^{1} \sqrt{1+\ln \Big(\phi  2^{\beta}
\Big[\frac{4\sqrt{2}\kappa}{s\|\mE\|_{2}}\Big]^{\beta/\gamma}\Big)} ds,
\end{align*}
\end{linenomath}
which, after some changes of variables, can be shown to be equal to
\begin{linenomath}
\begin{align*}
\Big(\frac{\beta}{\gamma}\Big)^{1/2}  
(\phi 2^{\beta}e)^{\gamma/\beta}
\Big(\frac{4\sqrt{2}\kappa}{\|\mE\|_{2}}\Big) 
\int_{\frac{\gamma}{\beta} [1 
+ \ln (\phi  2^{\beta}[\frac{4\sqrt{2}\kappa}
{\|\mE\|_{2}}]^{\beta/\gamma})]}^{\infty} 
\sqrt{x}e^{-x} dx.
\end{align*}
\end{linenomath}

Now, using the asymptotic behaviour of the (upper) incomplete Gamma function
\begin{linenomath}
\begin{align*}
\Gamma(y) = \int_{y}^{\infty} \sqrt{x}e^{-x} dx
\end{align*}
\end{linenomath}
as~$y\rightarrow\infty$, it is elementary to see that there exists
 a universal positive constant~$c_1$ such 
that~$\Gamma(y)\leq c_1\sqrt{y}e^{-y}$ for all~$y\geq 1/2$.
 Thus, since we are assuming~\eqref{cond_geq1}, we have
 \begin{linenomath}
\begin{align*}
J_{[\,]}\Big(1,\F,L_2\Big) \leq c_1 \sqrt{1 + \ln(\phi 2^{\beta}) 
+ \frac{\beta}{\gamma}\ln\Big(\frac{4\sqrt{2}\kappa}{\|\mE\|_{2}}\Big)}.
\end{align*}
\end{linenomath}

Finally, we use Theorem 2.14.2 of~\cite{VW96} to obtain that,
 if~\eqref{cond_geq1} holds, then (for a universal positive
 constant~$c_2$)
 \begin{linenomath}
\begin{align*}
E\Big[\sup_{x\in \Sigma}\Big|\sum_{k=1}^{n}\widehat{W}_k(x)\Big|
\Big] &\leq c_2 J_{[\,]}\Big(1,\F,L_2\Big) 
\|\mE\|_{2} \sqrt{n} \\
&\leq c_2 c_1 \sqrt{1 + \ln(\phi 2^{\beta}) 
+ \frac{\beta}{\gamma}\ln\Big(\frac{4\sqrt{2}\kappa}
{\|\mE\|_{2}}\Big)} \|\mE\|_{2} \sqrt{n},
\end{align*}
\end{linenomath}
which completes the proof with~$C_3=c_2 c_1$.
\end{proof}

We now investigate the~$L_2$-norm of the envelope~$\mE$. If we take, for example, an envelope of~$\F$ given by
\begin{linenomath}
\begin{align}
\mE(\omega) = 2\eps\sum_{k=\rho_1}^{\rho_2-1} \hat{\xi}_k , 
\mbox{ for any } \omega \in \mS,
\label{envelope}
\end{align}
\end{linenomath}
so that~$\mE\geq |f_x|$ for all~$f_x\in\F$ (recall the definition of an envelope function given at the beginning of Section~\ref{Est_F}), then we have

{\prop If~$\mE$ is given by~\eqref{envelope}, then it holds that
\begin{linenomath}
\begin{align*}
\|\mE\|_{2} = 4\sqrt{2}\eps.
\end{align*}
\end{linenomath}
\label{envelope_norm}
}

\begin{proof}
Use the fact that~$\sum_{k=\rho_1}^{\rho_2-1} \hat{\xi}_k$ is 
exponentially distributed with parameter~$\frac{1}{2}$.
\end{proof}

Thus, using Propositions~\ref{expect_R_tilde} and~\ref{envelope_norm}, we have under Assumptions~\ref{assump_1} and~\ref{assump_2}, for all~$n\in\IN$, 
\begin{linenomath}
\begin{align*}
E\Big[\sup_{x\in \Sigma}\Big|\sum_{k=1}^{n}\widehat{W}_k(x)\Big|
\Big] \leq F\sqrt{n}
\end{align*}
\end{linenomath}
where
\begin{linenomath}
\begin{align*}
F := C_4\eps \displaystyle\sqrt{1 + \ln(\phi 2^{\beta}) + \frac{\beta}{\gamma}\ln\Big(\frac{\kappa\vee (2\eps)}{\eps}\Big)},
\end{align*}
\end{linenomath}
and~$C_4$ is a universal positive constant.

\subsection{A tail bound involving~$\Psi$} 
\label{concentr}
 Recall the definition of~$\Psi$ from~\eqref{fct_psi}.
In this section, we will use a concentration inequality from~\cite{Adam08}, 
to estimate the tail of some useful random variable related to the numerator of $\sup_{x\in \Sigma}|\Psi(x)-1|$. 
Recalling the notation of Section \ref{Sim_MC}, we define, for~$x\in\Sigma$,
\begin{linenomath}
\begin{align*}
Z_n(x)=\sum_{i=1}^{n}\hat{\xi}_i-\hat{G}^X_n(x),
\end{align*}
\end{linenomath}
and we observe that $Z_n(\cdot)\eqlaw \sum_{i=1}^{|\mH^c|} \tilde{\xi}_i-\tilde{G}^{X}_{|\mH^c|}(\cdot)$.
Now we approximate~$Z_n$ by a sum of independent random elements. 
Specifically, we define
\begin{linenomath}
\begin{align*}
W_j(\cdot) = Z_{\rho_{j+1}-1}(\cdot) - Z_{\rho_j-1}(\cdot), \text{ for } j\geq 0,
\end{align*}
\end{linenomath}
with the convention that~$Z_0(\cdot)$=0, and we intend 
to approximate~$Z_n(\cdot)$ by a suitably chosen sum 
of $W_j(\cdot)$'s.

The random variables~$(\rho_{j+1} - \rho_j)_{j\geq 0}$ 
are i.i.d.\ Geometric($\frac{1}{2}$). The random 
elements~$(W_j(\cdot))_{j\geq 0}$ are independent and,
 additionally, the elements~$(W_j(\cdot))_{j\geq 1}$ 
are identically distributed. Also, for any~$x\in\Sigma$, we have
\begin{linenomath}
\begin{align*}
W_0(x) = \hat{\xi}_1(1 -\nu(x)) + \sum_{k=2}^{\rho_1-1} 
\hat{\xi}_k (1-I_k) (1 - \mu(\hat{x}_{k-1},x) ),
\end{align*}
\end{linenomath}
and, recalling~\eqref{widehatW},
\begin{linenomath}
\begin{align}
W_1(x) = \widehat{W}(x) = \sum_{k=\rho_1}^{\rho_2-1} \hat{\xi}_k (1-I_k) (1-\mu(\hat{x}_{k-1},x)).
\label{W1_law}
\end{align}
\end{linenomath}
Moreover, let us observe the fact that~$\IE[W_1(x)]=0$ for any~$x\in\Sigma$. To see this, first observe that
\begin{linenomath}
	\begin{align*}
    \IE\Big[ \sum_{k=\rho_1}^{\rho_2-1} \hat{\xi}_k (1-I_k) (1-\mu(\hat{x}_{k-1},x))\Big] = \IE\Big[\sum_{k=\rho_1+1}^{\rho_2-1} \big(1-\mu(\hat{x}_{k-1},x)\big)\Big],
	\end{align*}
\end{linenomath}
since~$I_{\rho_1}=1$, $I_k=0$ for~$k=\rho_1+1, \rho_1+2, \dots, \rho_2-1$, and the $(\hat{\xi}_k)_{k\geq 1}$ are Exponential($1$) distributed random variables, independent of all the other random elements. Then, by conditioning on the family~$(I_j)_{j\geq 1}$, we obtain
\begin{linenomath}
	\begin{align*}
	\IE\Big[\sum_{k=\rho_1+1}^{\rho_2-1} &\big(1-\mu(\hat{x}_{k-1},x)\big)\Big] \\
	&= \sum_{2\leq \ell < \ell'}\IE\Big[\ind_{\{\rho_1=\ell, \rho_2=\ell'\}}\sum_{k=\ell+1}^{\ell'-1}\IE\big[1-\mu(\hat{x}_{k-1},x) ~\big|~(I_j)_{j\geq 1}\big]\Big].
	\end{align*}
\end{linenomath}
But, given the regeneration times~$\rho_1$ and~$\rho_2$, $(\hat{x}_{k})_{k=\rho_1, \dots,\rho_2-2}$, is a Markov chain with starting law~$\Pi$ and transition density~$\mu$. Since $\Pi$ is also invariant for $\mu$ (that is, $\int \mu(x,y)\Pi(dy) = 1$, for all $x\in\Sigma$), we obtain that all the conditional expectations in the above display are null, so that~$\IE[W_1(x)]=0$.  

Let us denote, for~$n,m\in\IN$,
\begin{linenomath}
\begin{align}
\label{Rn}
R_n &= \sup_{x\in\Sigma}\Big|Z_n(x)\Big|, \nonumber\\
\tilde{R}_m &= \sup_{x\in\Sigma}\Big|\sum_{j=1}^m W_j(x)\Big|, 
\end{align}
\end{linenomath}
with the convention that~$\tilde{R}_0=0$.
Observe that in Section~\ref{Est_F} we actually proved that
\begin{linenomath}
\begin{align}
\IE[\tilde{R}_m] \leq F\sqrt{m},
\label{EPE_F_exp}
\end{align}
\end{linenomath}
where~$F$ is defined in~\eqref{F_exp}. 

We now obtain
{\prop    For all~$\theta>0$, it holds that
\begin{linenomath}
\begin{align*}
\IP\Big[R_n \geq 8 F \sqrt{n} + 7\theta n\Big] &  \leq 2\exp\Big\{-\frac{3}{128}\frac{\theta^2 n}{\eps^2}\Big\} + 6\exp\Big\{-3C_5\frac{\theta n}{\eps\ln(3n+1)}\Big\}\\
&\phantom{***} +  2\exp\Big\{-\frac{\theta n}{8\eps}\Big\},
\end{align*}
\end{linenomath}
where $C_5$ is a universal positive constant.
\label{Prop_Rn}
}

\begin{proof}
We will use an argument analogous to the one used in the proof of Lemma 2.9 in \cite{CP}. To begin, let us assume that there exists a positive constant~$C_5$ such that, for all~$\theta>0$,
\begin{linenomath}
\begin{align}
\IP\Big[\tilde{R}_m \geq \frac{3}{2} F \sqrt{m} + \theta m\Big] \leq \exp\Big\{-\frac{\theta^2 m}{128\eps^2}\Big\} 
+ 3\exp\Big\{-C_5\frac{\theta m}{\eps\ln(m+1)}\Big\}.
\label{assump}
\end{align}
\end{linenomath}
(This statement will be proved later, in Proposition~\ref{Prop_Rn_tilde}).

Using Markov's inequality and~\eqref{EPE_F_exp}, it is elementary to show that, for~$c_1= 4/3$,
\begin{linenomath} 
\begin{align*}
\IP\Big[\tilde{R}_m \geq \frac{3}{2} c_1 F \sqrt{m} + \theta m\Big] \leq \frac{1}{2}, 
\end{align*}
\end{linenomath}
so that 
\begin{linenomath}
\begin{align}
\IP\Big[\tilde{R}_m \geq& 2 F \sqrt{m} + \theta m\Big] \nonumber \\ 
&\leq \frac{1}{2} \wedge \Big(\exp\Big\{-\frac{\theta^2 m}{128\eps^2}\Big\} 
+ 3\exp\Big\{-C_5\frac{\theta m}{\eps\ln(m+1)}\Big\}\Big).
\label{EP2} 
\end{align}
\end{linenomath}

Now, let us define 
\begin{linenomath}
\begin{align*}
\tilde{M} = \min\Big\{ i\in[0,3m] : \tilde{R}_i \geq 8 F \sqrt{m} + 6\theta m\Big\}, 
\end{align*}
\end{linenomath}
with the convention that $\min\emptyset=\infty$, so that 
\begin{linenomath}
\begin{align*}
\IP\Big[\tilde{M}\in[0,3m]\Big] = \IP\Big[\max_{i\in[0,3m]}\tilde{R}_i \geq& 8 F \sqrt{m} + 6\theta m\Big].
\end{align*}
\end{linenomath}
Then, using~\eqref{EP2} one gets
\begin{linenomath}
\begin{align*}
\exp\Big\{&-\frac{3}{128}\frac{\theta^2 m}{\eps^2}\Big\} + 3\exp\Big\{-3C_5\frac{\theta m}{\eps\ln(3m+1)}\Big\} \\
&\geq \IP\big[\tilde{R}_{3m} \geq 2 F \sqrt{3m} + 3\theta m\big] \\
&= \sum_{j=0}^{3m} \IP[\tilde{M}=j] \IP\big[\tilde{R}_{3m} \geq 2 F \sqrt{3m} + 3\theta m ~\big|~ \tilde{M}=j\big] \\
&\geq \sum_{j=0}^{3m} \IP[\tilde{M}=j] \IP\big[\tilde{R}_{3m-j} < 2 F \sqrt{3m} + 3\theta m\big] \\
&\geq \IP\big[\tilde{M}\in[0,3m]\big] \min_{j\in[0,3m]} \IP\Big[\tilde{R}_{3m-j} < 2 F \sqrt{3m-j} + \theta (3m-j)\Big] \\
&\geq \frac{1}{2} \IP\big[\tilde{M}\in[0,3m]\big],
\end{align*}
\end{linenomath}
where, to obtain the second inequality, we use the Markov property of the random walk~$(\sum_{j=1}^m W_j)_m$ and the fact that, to be above the level~$2 F \sqrt{3m} + 3\theta m$ at time~$3m$, being above the level~$8 F \sqrt{m} + 6\theta m$ at time $j$, it suffices that the process~$(\tilde{R}_{n})$ decreases less than~$2 F \sqrt{3m} + 3\theta m$ during the time interval $3m-j$. This in turn proves that
\begin{linenomath}
\begin{align}
\IP\Big[\max_{i\in[0,3m]} & \tilde{R}_i \geq 8 F \sqrt{m} + 6\theta m\Big] \nonumber \\
&\leq 2\exp\Big\{-\frac{3}{128}\frac{\theta^2 m}{\eps^2}\Big\} + 6\exp\Big\{-3C_5\frac{\theta m}{\eps\ln(3m+1)}\Big\}. 
\label{Max_inqty}
\end{align}
\end{linenomath}

Now let us define
\begin{linenomath}
\begin{align*}
\sigma_n = \min\{j\geq 1 : \rho_j > n\},
\end{align*}
\end{linenomath}
so that~$\sigma_n - 1$ represents the number of regenerations of the Markov chain~$X$ until time~$n$. Observe that~$\sigma_n - 1$ is a Binomial($n-1,\frac{1}{2}$) distributed random variable. Then, under Assumption~\ref{assump_3}, using the triangular inequality and~\eqref{Rn}, we obtain 
\begin{linenomath}
\begin{align}
R_n \leq \tilde{R}_{\sigma_n-1} + 2\eps\sum_{i=1}^{\rho_1-1}\hat{\xi}_i + 2\eps\sum_{i=n+1}^{\rho_{\sigma_n}} \hat{\xi}_i.
\label{Rn_inqty}
\end{align}
\end{linenomath}
The last two terms in the above inequality take care of the terms~$W_0$ and $W_{\sigma_n-1}$, respectively, when comparing~$R_n$ and~$\tilde{R}_{\sigma_n-1}$.  
Note that~$\rho_{\sigma_n}-n$ and~$\rho_1-1$ are both Geometric($\frac{1}{2}$) distributed random variables. Consequently, $\sum_{i=n+1}^{\rho_{\sigma_n}} \hat{\xi}_i$ and~$\sum_{i=1}^{\rho_1-1} \hat{\xi}_i$ are both Exponential($\frac{1}{2}$) distributed and we deduce that
\begin{linenomath}
\begin{align}
\IP\Big[ 2\eps\sum_{i=n+1}^{\rho_{\sigma_n}} \hat{\xi}_i + 2\eps\sum_{i=1}^{\rho_1-1} \hat{\xi}_i \geq \theta n\Big] 
&\leq \IP\Big[ 2\eps\sum_{i=n+1}^{\rho_{\sigma_n}}\hat{\xi}_i \geq \frac{\theta n}{2}\Big] + \IP\Big[2\eps\sum_{i=1}^{\rho_1-1} \hat{\xi}_i \geq \frac{\theta n}{2}\Big] \nonumber\\
&= 2\exp\Big\{-\frac{\theta n}{8\eps}\Big\}.
\label{Exp_inqty}
\end{align}
\end{linenomath}

Finally, using~\eqref{Rn_inqty} together 
with~\eqref{Max_inqty} and~\eqref{Exp_inqty}, we obtain that
\begin{linenomath}
\begin{align*}
\IP\big[& R_n \geq 8 F \sqrt{n} + 7\theta n\big] \\
& \leq \IP\Big[ \max_{i\in[0,3n]}\tilde{R}_i \geq 8 F \sqrt{n}
 + 6\theta n\Big]  + \IP\Big[ \sum_{i=n+1}^{\rho_{\sigma_n}} \hat{\xi}_i 
+ \sum_{i=1}^{\rho_1-1} \hat{\xi}_i \geq \frac{\theta n}{2\eps}\Big] \\
&\leq 2\exp\Big\{-\frac{3}{128}\frac{\theta^2 n}{\eps^2}\Big\}
 + 6\exp\Big\{-3C_5\frac{\theta n}{\eps\ln(3n+1)}\Big\} 
+  2\exp\Big\{-\frac{\theta n}{8\eps}\Big\} .
\end{align*}
\end{linenomath}
 \end{proof}
 
 Then, we deduce the following
 
 {\cor   For all $i,n\in\IN$, it holds that
 \begin{linenomath}
\begin{align*}
 \IP\big[R_n \geq 8 F \sqrt{n} (1+i)\big] &  
\leq 2\exp\Big\{-C_6\frac{F^2 i^2 }{\eps^2}\Big\} 
+ 8\exp\Big\{-C_7\frac{Fi}{\eps}\frac{\sqrt{n}}{\ln(3n+1)}\Big\},
\end{align*}
\end{linenomath}
where $C_6$ and~$C_7$ are universal positive constants.
\label{Cor_Rn}
 }
 
 \begin{proof}
Take~$\theta=\frac{8}{7}\frac{Fi}{\sqrt{n}}$ in 
Proposition~\ref{Prop_Rn}.
 \end{proof}

Before proving assertion~\eqref{assump},
which was assumed to be true in the beginning of the proof 
of Proposition~\ref{Prop_Rn}, we must prove some 
preliminary results.

{\lem It holds that
\begin{linenomath}
\begin{align*}
\sigma^2 := \sup_{x\in\Sigma} \sum_{j=1}^m \IE\big[W_j(x)^2\big] 
\leq 32\eps^2m.
\end{align*}
\end{linenomath}
\label{lem_sigma2}
}

\begin{proof}
Since, for any~$x\in\Sigma$, the elements of~$(W_j(x))_{j\geq 1}$
 are i.i.d., we only have to prove that, for any~$x\in\Sigma$,
 \begin{linenomath}
\begin{align*}
\IE\big[W_1(x)^2\big] \leq 32\eps^2.
\end{align*}
\end{linenomath}

We use~\eqref{W1_law} and~\eqref{max_eps} to write, for any~$x\in\Sigma$,
\begin{linenomath}
\begin{align*}
\IE\big[W_1(x)^2\big] \leq  (2\eps)^2 
\IE\Big[\Big(\sum_{k=\rho_1}^{\rho_2-1} \hat{\xi}_k \Big)^2 \Big].
\end{align*} 
\end{linenomath}
Since~$\rho_2-\rho_1$ is Geometric($\frac{1}{2}$), we use the fact that~$\sum_{k=\rho_1}^{\rho_2-1} \hat{\xi}_k$ is exponentially 
distributed with parameter~$\frac{1}{2}$ to see that 
\begin{linenomath}
\begin{align*}
\IE\big[W_1(x)^2\big] \leq  32\eps^2,
\end{align*}
\end{linenomath}
and this completes the proof.
\end{proof}

In order to formulate the next lemma, 
we define the so-called $\psi_1$-Orlicz norm of a 
random variable~$X$, in the following way:
\begin{linenomath}
\begin{align*}
\|X\|_{\psi_{1}} = \inf\big\{t>0 : \IE e^{|X|/t} \leq 2\big\},
\end{align*}
\end{linenomath}
see Definition 1 of~\cite{Adam08}.

{\lem It holds that
\begin{linenomath}
\begin{align*}
\Big\| \max_{1\leq j\leq m} \sup_{x\in\Sigma}|W_j(x)| 
\Big\|_{\psi_1} \leq C_8 \eps  \ln(m+1),
\end{align*}
\end{linenomath}
where~$C_8$ is a universal positive constant.
\label{lem_orlicz}
}

\begin{proof}
First, Lemma 2.2.2 of~\cite{VW96} provides the inequality
\begin{linenomath}
\begin{align*}
\Big\|\max_{1\leq j \leq m} \sup_{x\in\Sigma} |W_j(x)| 
\Big\|_{\psi_{1}} \leq c_1 \max_{1\leq j \leq m} \Big\|\sup_{x\in\Sigma} 
|W_j(x)| \Big\|_{\psi_{1}} \ln (m+1),
\end{align*}
\end{linenomath}
for a universal positive constant~$c_1$.

But, due to~\eqref{max_eps},
\begin{linenomath}
\begin{align*}
\sup_{x\in\Sigma}\Big|\sum_{k=\rho_1}^{\rho_2-1} \hat{\xi}_k (1-I_k) (1-\mu(\hat{x}_{k-1},x)) \Big| \leq 2\eps \sum_{k=\rho_1}^{\rho_2-1} \hat{\xi}_k,
\end{align*}
\end{linenomath}
so that (recall~\eqref{W1_law})
\begin{linenomath}
\begin{align*}
\Big\|\sup_{x\in\Sigma} |W_1(x)| \Big\|_{\psi_1} \leq 2\eps 
\Big\|\sum_{k=\rho_1}^{\rho_2-1} \hat{\xi}_k \Big\|_{\psi_1}.
\end{align*}
\end{linenomath}
Then, since~$\sum_{k=\rho_1}^{\rho_2-1} \hat{\xi}_k$ is exponentially 
distributed with mean~$2$ and the~$\psi_1$-Orlicz norm of an 
exponential random variable equals twice its mean, we
obtain the result (recalling that the~$W_j$'s, for $j\geq 1$, 
are independent and identically distributed).
\end{proof}

Now, in order to address the problem of estimating the probability 
involving~$\tilde{R}_m$ (whose bound was postulated in~\eqref{assump})
from the viewpoint of the theory of empirical processes, we recall 
the space~$\mS=\R^{\Sigma}$ and the class~$\F=(f_x)_{x\in\Sigma}$ of 
functions~$f_x:\mS\rightarrow \R$ such that~$f_x(\omega)=\omega(x)$, 
for~$\omega\in\mS$, so that the above mentioned probability can be 
rewritten as
\begin{linenomath}
\begin{align}
\IP\Big[ \sup_{f_x\in\F} \Big|\sum_{j=1}^{m} f_x(W_j(\cdot))\Big| 
\geq \frac{3}{2} F \sqrt{m} + \theta m \Big],
\label{prob_int}
\end{align}
\end{linenomath}
with~$W_j(\cdot)$ being interpreted as a vector in~$\mS$ whose components 
are~$W_j(x)$ for each~$x\in\Sigma$. In this setting, we are able to apply 
Theorem~4 of~\cite{Adam08} to prove~\eqref{assump}, 
and this is done in the next proposition.

{\prop There exists a universal positive constant~$C_5$ such that, 
for all~$\theta>0$,
\begin{linenomath}
\begin{align*}
\IP\Big[\tilde{R}_m \geq \frac{3}{2} F \sqrt{m} + \theta m\Big] 
\leq \exp\Big\{-\frac{\theta^2 m}{128\eps^2}\Big\} 
+ 3\exp\Big\{-C_5\frac{\theta m}{\eps\ln(m+1)}\Big\}.
\end{align*}
\end{linenomath}
\label{Prop_Rn_tilde}
}

\begin{proof}
We use~\eqref{EPE_F_exp} and just apply Theorem 4 of~\cite{Adam08} 
(with~$\delta=1$, $\eta=1/2$ and~$\alpha=1$ there) 
together with Lemmas~\ref{lem_sigma2} and~\ref{lem_orlicz}, 
to see that there exist universal positive constants~$C$ and~$C_8$
($C$ is from Theorem 4 of~\cite{Adam08} and~$C_8$ is 
from Lemma~\ref{lem_orlicz}) such that, for all~$t>0$,
\begin{linenomath}
\begin{align*}
\IP\Big[\tilde{R}_m \geq \frac{3}{2} F \sqrt{m} + t \Big] 
\leq \exp\Big\{-\frac{t^2}{128\eps^2 m}\Big\} 
+ 3\exp\Big\{-\frac{t}{C_8C\eps\ln(m+1)}\Big\}.
\end{align*}
\end{linenomath}
We conclude the proof by setting~$t=\theta m$, 
for~$\theta>0$, and~$C_5=(C_8C)^{-1}$.
\end{proof}

\subsection{Proof of Proposition \ref{goodenv}} \label{Proof_Av}

We begin this section obtaining a tail estimate for the cardinality of the random set~$\mH$ introduced in~\eqref{set_H}, which verifies
\begin{linenomath}
\begin{align}
\label{Xidef}
|\mH| = 
\sum_{j=2}^{n-1} (I_j I_{j+1})  + I_n.
\end{align}
\end{linenomath}

{\prop \label{Prop_H_comp} There exist a positive universal
constant~$C_9$ such that, for all $n\geq 44$,
 it holds that
 \begin{linenomath}
\begin{align*}
\IP[\C^c] := \IP\Big[|\mH|\leq \frac{n}{24}\Big] \leq e^{-C_9 n} \leq \frac{1}{2}.
\end{align*}
\end{linenomath}
}
\begin{proof}
First, observe that for $n\geq 4 $ we have that
\begin{linenomath}
\begin{align*}
\IP\Big[|\mH| \leq \frac{n}{24} \Big] 
&\leq \IP\Big[\sum_{j=2}^{n-1} (I_j I_{j+1})  
\leq \frac{n}{24} \Big] 
\leq \IP\Big[\sum_{k=1}^{\lfloor\frac{n-2}{2}\rfloor}
 (I_{2k} I_{2k+1})  \leq \frac{n}{24} \Big].
\end{align*}
\end{linenomath}
Moreover, observe that the random variables~$(I_{2k} I_{2k+1})$, 
$k=1,2,\dots,\lfloor\frac{n-2}{2}\rfloor$, are independent 
 Bernoulli($\frac{1}{4}$). Now, we recall the standard lower tail bound
 for the binomial law: for~$\mathcal{X}\sim$ Binomial($m,p$) 
and~$\delta\geq 0$, we have that
\begin{linenomath}
\begin{equation*}
\IP\big[\mathcal{X}\leq (1-\delta) m p\big]  \leq \exp\{-m \II(p,\delta)\},
\end{equation*}
\end{linenomath}
where 
\begin{linenomath}
\begin{equation*}
\II(p,\delta):=p(1-\delta)\ln(1-\delta)
+p\Big(\frac{1-p}{p}+\delta\Big)
\ln\Big(1+\frac{\delta p}{1-p}\Big).
\end{equation*}
\end{linenomath}
Applying the above formula to the random variable
\begin{linenomath}
\begin{align*}
\sum_{k=1}^{\lfloor\frac{n-2}{2}\rfloor} (I_{2k} I_{2k+1})
\end{align*}
\end{linenomath}
with~$\delta=1/2$, we obtain that
\begin{linenomath}
\begin{align*}
\IP\Big[\sum_{k=1}^{\lfloor\frac{n-2}{2}\rfloor} (I_{2k} I_{2k+1}) 
 \leq   \frac{n}{24} \Big] \leq \exp\Big\{-\frac{n}{3}\II\Big(\frac{1}{4},\frac{1}{2}\Big)\Big\}\leq \frac{1}{2},
\end{align*}
\end{linenomath}
for~$n\geq 44$.
\end{proof}


Next, we obtain an upper bound for 
\[
 \IE\Big[1\wedge\sup_{x\in \Sigma}|\Psi(x)-1|^3 \; \Big|\; \C ~\Big].
\]
Then, the upper bound for $\IP[\A_i^c\mid \C]$ 
will be a direct application of Markov's inequality.
First observe that 
\begin{linenomath}
\begin{align}
\label{RRR}
 \IE\Big[1\wedge\sup_{x\in \Sigma}|\Psi(x)-1|^3 \; \Big|\; \C ~\Big]&\leq  \IE\Big[\big(1\wedge\sup_{x\in \Sigma}|\Psi(x)-1|^3\big)\ind_{\BB} \; \Big|\; \C ~\Big]+\IP[\BB^c\mid \C]
\end{align}
\end{linenomath}
where 
\begin{equation}
\BB=\Big\{\sup_{x,y\in \Sigma}|\tilde{G}^X_{|\mH^c|}(x)-\tilde{G}^X_{|\mH^c|}(y)|\leq \frac{\sum_{i=|\mH^c|+1}^{n}\tilde{\xi}_i}{2}\Big\}.
\end{equation}

Applying Proposition~\ref{Prop_H_comp} and recalling~\eqref{Rn}, we have
\begin{linenomath}
\begin{align}
\label{EXPPSI}
\IE\Big[\big(1\wedge \sup_{x\in \Sigma}|\Psi(x)-1|^3\big)\ind_\BB\; \Big|\; \C~\Big]
&\leq 2^7\IE\Bigg[\sup_{x,y\in \Sigma}\Big|
\frac{\tilde{G}^X_{|\mH^c|}(x)-\tilde{G}^X_{|\mH^c|}(y)}
{\sum_{i=|\mH^c|+1}^{n}\tilde{\xi}_i}\Big|^3 
\ind_{\C}\Bigg]\nonumber\\
&\leq 2^7\IE\Big[\sup_{x,y\in\Sigma}|
\tilde{G}^X_{|\mH^c|}(x)-\tilde{G}^X_{|\mH^c|}(y)|^3\Big]
\IE\Big[\Big(\sum_{i=1}^{\lfloor
 \frac{n}{24}\rfloor}\tilde{\xi}_i\Big)^{-3}\Big]\nonumber\\
&= 2^7\IE\Big[\sup_{x,y\in\Sigma}|
Z_n(x)-Z_n(y)|^3\Big]
\IE\Big[\Big(\sum_{i=1}^{\lfloor
	\frac{n}{24}\rfloor}\tilde{\xi}_i\Big)^{-3}\Big]\nonumber\\
&\leq 2^{10}\IE\Big[\sup_{x\in\Sigma}|
Z_n(x)|^3\Big]
\IE\Big[\Big(\sum_{i=1}^{\lfloor
	\frac{n}{24}\rfloor}\tilde{\xi}_i\Big)^{-3}\Big]\nonumber\\
&=2^{10}\IE[R_n^3]\IE\Big[\Big(
\sum_{i=1}^{\lfloor  \frac{n}{24}\rfloor}\tilde{\xi}_i\Big)^{-3}\Big],
\end{align}
\end{linenomath}
where in the second step, we used the fact that, conditionally 
on $\sigma(I_j,j\geq 1)$ the numerator and the denominator of
 the ratio in the first line are independent, and the event
 $\C = \{|\mH|>  \frac{n}{24}\}$ is measurable with respect to 
$\sigma(I_j,j\geq 1)$. In the third step, we use the fact that $Z_n(\cdot)\eqlaw \sum_{i=1}^{|\mH^c|} \tilde{\xi}_i-\tilde{G}^{X}_{|\mH^c|}(\cdot)$.

Then, using an integration by parts,
 Corollary~\ref{Cor_Rn}, and the fact that the square root 
in~(\ref{F_exp}) is greater than one,  we obtain that
\begin{linenomath}
\begin{align*}
\IE\Big[\Big(\frac{R_n}{8F\sqrt{n}}\Big)^3\Big]&=3\int_0^{\infty}t^2
\IP\Big[\frac{R_n}{8F\sqrt{n}}>t\Big]dt\nonumber\\
&\leq 1+3\int_1^{\infty}t^2
\IP\Big[\frac{R_n}{8F\sqrt{n}}>t\Big]dt\nonumber\\
&\leq c_1,
\end{align*}
\end{linenomath}
where $c_1$ is positive. Therefore, we have 
\begin{equation}
\label{RICO}
\IE[R_n^3]\leq c_2F^3n^{3/2},
\end{equation}
where $c_2$ is positive.
 On the other hand, since 
$\big(\sum_{i=1}^{\lfloor \frac{n}{24}\rfloor}\tilde{\xi}_i\big)^{-1}$ 
is an Inverse Gamma random variable with parameters 
$(\lfloor \frac{n}{24}\rfloor,1)$, we obtain that
\begin{equation}
\label{InvGam}
\IE\Big[\Big(\sum_{i=1}^{\lfloor \frac{n}{24}\rfloor}\tilde{\xi}_i\Big)^{-3}\Big]= \frac{1}{(\lfloor \frac{n}{24}\rfloor-1)(\lfloor \frac{n}{24}\rfloor-2)(\lfloor \frac{n}{24}\rfloor-3)}\leq c_3 n^{-3}
\end{equation}
for $n\geq 96$ and $c_3>0$.
Gathering  (\ref{EXPPSI}), (\ref{RICO}), and~(\ref{InvGam}) 
we obtain, for~$n\geq 96$,
\begin{linenomath}
\begin{equation}
\label{RRR1}
\IE\Big[\big(1\wedge \sup_{x\in \Sigma}|\Psi(x)-1|^3\big)\ind_\BB\; \Big|\; \C~\Big]\leq c_4F^3n^{-3/2}
\end{equation}
\end{linenomath}
for some positive constant $c_4$. 

The term $\IP[\BB^c\mid \C]$ of (\ref{RRR}) can be treated using the Markov inequality and the above estimates to obtain for $n\geq 96$,
\begin{linenomath}
	\begin{equation}
	\label{RRR2}
\IP[\BB^c\mid \C]\leq c_5 F^3n^{-3/2}
	\end{equation}
\end{linenomath}
for some positive constant $c_5$. 

Finally, gathering (\ref{RRR}), (\ref{RRR1}), (\ref{RRR2}) and applying the Markov inequality, we deduce that for $n\geq 96$,
\begin{equation*}
\IP[\A_i^c\mid \C]\leq \frac{C_2}{(i+1)^3}
\end{equation*}
for some positive constant $C_2$.
\qed

\section{Proof of Theorem~\ref{Main_Thm}}
 \label{Main_Thm_proof}
We estimate $\IP[\Upsilon^c]$ from above
 (recall that $\Upsilon$ is the coupling event from 
Section~\ref{coupling}) to obtain an upper bound on the total 
variation distance between $L_n^X$ and $L_n^Y$.
At this point, we mention that we will use the notation from Section~\ref{coupling}.
By definition of the total variation distance, we have that
\begin{equation}
\label{Theoproof1}
\text{d}_{\text{TV}}(L_n^X , L_n^Y)\leq \IP[\Upsilon^c].
\end{equation}

First, let us decompose $\Upsilon$ according to
 $\mathsf{C}=\{|\mH|> \frac{n}{24}\}$ (recall~\eqref{event_C}) 
and its complement:
\begin{equation}
\label{Theoproof2}
\IP[\Upsilon^c]=\IP[\Upsilon^c, \mathsf{C}]+\IP[\Upsilon^c, \mathsf{C}^c].
\end{equation}
Now, we partition $\C$ using the events $\BB_1:=\A_1$ and $\BB_{i+1}:=\A_{i+1}\setminus\A_i$, for $i\geq 1$, to write
\begin{linenomath}
\begin{equation*}
\IP[\Upsilon^c,\mathsf{C}]=\sum_{i=1}^{\infty}\IP[\Upsilon^c,\mathsf{B}_i].
\end{equation*}
\end{linenomath}
Since $\mathsf{B}_i$ is $\sigma(\W)$-measurable for any~$i\geq 1$ (we recall that $\W$ was introduced in 
Section~\ref{coupling}), we have that
\begin{linenomath}
\begin{equation*}
\IP[\Upsilon^c, \mathsf{B}_i]=\IE[{\bf 1}_{\mathsf{B}_i}\IP[\Upsilon^c\mid \W]].
\end{equation*}
\end{linenomath}
Then, observe that from the coupling construction
of Section~\ref{coupling}, we have 
\begin{linenomath}
\begin{equation*}
\IP[\Upsilon^c\mid \W]\leq \big\|\IP[{V}\in \cdot\mid \W]-\IP[{V'}\in \cdot\mid \W]\big\|_{\text{TV}}.
\end{equation*}
\end{linenomath}
Hence, applying Proposition \ref{Propmulti} to the term in the 
right-hand side with $\delta_0=1$ we obtain, on the sets
 $\mathsf{B}_i$,
 \begin{linenomath}
\begin{equation*}
\IP[\Upsilon^c\mid \W]\leq C_1(1)(1+i)F.
\end{equation*}
\end{linenomath}
Recalling~\eqref{event_Av}, we deduce that
\begin{linenomath}
\begin{equation*}
\IP[\Upsilon^c, \mathsf{C}]\leq C_1(1)F\sum_{i=1}^{\infty}(i+1)\IP[\mathsf{B}_i]\leq C_1(1)F\Big(2+\sum_{i=2}^{\infty}(i+1)\IP[\mathsf{A}^c_{i-1}\mid \mathsf{C}]\Big),
\end{equation*}
\end{linenomath}
 which by Proposition \ref{goodenv} implies, for $n\geq 96$,

\begin{equation}
\label{Theoproof3}
\IP[\Upsilon^c, \mathsf{C}]\leq c_1 F
\end{equation}
for some positive constant~$c_1$. 

Regarding the second term of the sum in~\eqref{Theoproof2}, we can write
\begin{linenomath}
\begin{align*}
\IP[\Upsilon^c, \mathsf{C}^c] = \IE[\ind_{\C^c}\IP[\Upsilon^c\mid \I]],
\end{align*}
\end{linenomath}
since~$\C^c$ is $\sigma(\I)$-measurable (recall that~$\I=(I_1,\dots,I_n)$). By construction of our coupling, observe that $\IP[\Upsilon^c\mid \I]\leq 2\varepsilon n$, on~$\C^c$. Thus, we obtain
\begin{linenomath}
\begin{align*}
\IP[\Upsilon^c, \mathsf{C}^c] \leq 2\varepsilon n \IP[\C^c].
\end{align*}
\end{linenomath}
Using Proposition~\ref{Prop_H_comp} we have that
\begin{linenomath}
\begin{align}
\label{expon_bound}
\IP[\Upsilon^c, \mathsf{C}^c] \leq 2\varepsilon n e^{-C_9 n},
\end{align}
\end{linenomath}
for~$n\geq 44$.

Finally, combining~(\ref{Theoproof1}), 
(\ref{Theoproof2}), 
(\ref{Theoproof3}), (\ref{expon_bound}) 
and~(\ref{F_exp}), 
we obtain Theorem~\ref{Main_Thm} for $n\geq 96$. 

For $n<96$, we simply perform a step-by-step coupling between 
the Markov chain~$X$ and the sequence~$Y$ 
(as described in the introduction) 
to obtain $\text{d}_{\text{TV}}(L_n^X,L_n^Y)\leq 95 \eps$ 
and thus prove Theorem~\ref{Main_Thm}.
\qed

\section{Proof of Theorem~\ref{Thm2}}
 \label{Second_Thm}

First, we prove a preliminary lemma. As in Section~\ref{TV_binomial},
 consider again two binomial point processes on some measurable space
 $(\Omega, \mathcal{T})$ with laws~${\P}_n$ and~${\Q}_n$ 
of respective parameters $({\p}_n,n)$ and $(\q_n,n)$, 
where $n\in \mathbb{N}$ and~$\p_n$, $\q_n$ are 
two probability laws on $(\Omega, \mathcal{T})$ such that
 $\q_n\ll \p_n$. Then, we have 

{\lem 
\label{Propmulti_2}
Let~$\delta>0$ such that, for all $n\in \mathbb{N}$, 
$|\frac{\mathrm{d}\q_n}{\mathrm{d}\p_n}(x)-1|\leq\delta n^{-1/2}$ 
for all~$x\in \Omega$. Then
\begin{linenomath}
\begin{equation*}
\sup_{n\geq 1}\|\P_n-\Q_n\|_{\emph{TV}}\leq 1 - C_{10}(\delta),
\end{equation*}
\end{linenomath}
where~$C_{10}(\delta)$ is a positive constant depending on $\delta$.
}

\begin{proof}
As pointed out in the proof of Proposition~\ref{Propmulti}, 
$\P_n$ and~$\Q_n$ can be seen as probability measures 
on the space of $n$-point measures
\[
 \M_n=\big\{m:m=\sum_{i=1}^{n}\bdelta_{x_i},
    x_i\in \Omega, 1\leq i\leq n\big\}
\]
 endowed with the $\sigma$-algebra generated by the mappings 
$\Phi_B:\M_n\to \Z_+$ defined by $\Phi_B(m)=m(B)
 =\sum_{i=1}^n\bdelta_{x_i}(B)$, 
for all $B\in \mathcal{T}$. Also, recall that $\Q_n\ll \P_n$ 
and its Radon-Nikodym derivative with respect to $\P_n$ is given by
\[
\frac{\text{d}\Q_n}{\text{d}\P_n}(m)
=\prod_{i=1}^n\frac{\text{d}\q_n}{\text{d}\p_n}(x_i)
\]
where $m=\sum_{i=1}^{n}\bdelta_{x_i}$. 
Moreover, for~$n\in \mathbb{N}$, recall 
the functions~$f_n,g_n:\Omega\to \R$ given by
\begin{linenomath}
\begin{align*}
f_n(x) = \frac{\text{d}\q_n}{\text{d}\p_n}(x)-1 \ \ \mbox{ and } \ \ g_n(x) = \ln(f_n(x)+1), \ \ \mbox{ for } x\in \Omega.
\end{align*}
\end{linenomath}
We start  by proving the lemma for all large enough~$n$.

It is convenient to introduce now two new distinct elements~$\0_1$ and~$\0_2$ in order to define a new space~$\hat{\Omega}=\Omega\cup\{\0_1,\0_2\}$ (we assume that 
$\0_1, \0_2\notin \Omega$), endowed with the
$\sigma$-algebra~$\hat{\mathcal{T}}:=\sigma(\mathcal{T},\{\0_1\})$. 
Then, on~$(\hat{\Omega},\hat{\mathcal{T}})$ we consider a 
 new binomial point process with law~${\hP}_{n,k}$ of parameters
 $(\hp_n,k)$, where $k\in \mathbb{N}$ and~$\hp_n$
 is the probability law on $(\hat{\Omega}, \hat{\mathcal{T}})$ given by
 \begin{linenomath}
\begin{align*}
\hp_n(A) = 
\begin{cases}
\frac{\p_n(A)}{2}, &\mbox{ for } A\in\mathcal{T}, \\
\frac{1}{4}, &\mbox{ for } A\in \{\{\0_1\},\{\0_2\}\},
\end{cases}
\end{align*}
\end{linenomath}
Additionally, for $n>\delta^2 $, consider another binomial point 
process on~$(\hat{\Omega},\hat{\mathcal{T}})$ with
 law~${\hat{\Q}}_{n,k}$ and parameters~$(\hq_n,k)$, 
where~$\hq_n$ is the probability law on $(\hat{\Omega}, \hat{\mathcal{T}})$ such that~$\hq_n(A)=\frac{\q_n(A)}{2}$ 
for all~$A\in\mathcal{T}$ and
\begin{linenomath}
\begin{align*}
\hq_n(\{\0_1\})=\frac{1}{4}\Big(1+\frac{\delta}{\sqrt{n}}\Big),\; \hq_n(\{\0_2\}) = \frac{1}{4}\Big(1-\frac{\delta}{\sqrt{n}}\Big)
\end{align*}
\end{linenomath}
so that $\hq_n(\{\0_1\})+\hq_n(\{\0_2\}) = 1/2$.
Thus, $\hP_{n,k}$ and~$\hQ_{n,k}$ can be seen as probability 
measures on the space
\[
\hat{\M}_k=\big\{\hat{m}:\hat{m}=\sum_{i=1}^{k}\bdelta_{x_i}, 
x_i\in \hat{\Omega}, 1\leq i\leq k\big\}.
\]

We need to introduce the corresponding functions~$\hf_n,\hg_n,:\hat{\Omega}\to \R$ given by
\begin{linenomath}
\begin{align*}
\hf_n(x) = \frac{\text{d}\hq_n}{\text{d}\hp_n}(x)-1 \ \ 
\mbox{ and } \ \ \hg_n(x) = \ln(\hf_n(x)+1).
\end{align*}
\end{linenomath}
Also, define~$\hh_n:\hat{\Omega}\to \R$ as~$\hh_n = \delta^{-1}n^{1/2}\hf_n$ and the set~$\hat{\M}^{\hat{\Q}}_k
=\big\{\hat{m}\in\hat{\M}_k:\frac{\text{d}\hat{\Q}_{n,k}}
{\text{d}\hP_{n,k}}(\hat{m}) \geq 1\big\}$.
 Since
 \begin{linenomath}
\begin{align*}
\frac{\text{d}\hat{\Q}_{n,k}}{\text{d}\hP_{n,k}}(\hat{m}) 
\geq 1 \Leftrightarrow \sum_{i=1}^k 
 \ln\Big(\frac{\text{d}\hat{\q}_n}{\text{d}\hp_n}(x_i)\Big) 
\geq 0, 
\end{align*}
\end{linenomath}
for~$\hat{m}\in\hat{\M}_k$, we have for any~$n,k\in\IN$
\begin{linenomath}
\begin{align*}
\|\hP_{n,k}-\hat{\Q}_{n,k}\|_{\text{TV}}
&=\int_{\hat{\M}_k}\Big(\frac{\text{d}\hat{\Q}_{n,k}}
{\text{d}\hP_{n,k}}(\hat{m})-1\Big)^+\text{d}\hP_{n,k}(\hat{m})\nonumber\\
&=\int_{\hat{\M}^{\hat{\Q}}_k}\Big(\frac{\text{d}\hat{\Q}_{n,k}}
{\text{d}\hP_{n,k}}(\hat{m})-1\Big)\text{d}\hP_{n,k}(\hat{m})\nonumber\\
&\leq 1 - \hP_{n,k}\Big[\hat{m}(\hat{g}_n) \geq 0\Big].
\end{align*}
\end{linenomath}
Next, we will bound~$\hP_{n,2n}[\hat{m}(\hat{g}_n) \geq 0]$ from above.

If we define~$n_1=n_1(\delta)=\lceil 3\delta^2\rceil$ then using the 
fact that~$\ln(1+x)\geq x-x^2$ for~$x\in(-1/\sqrt{3},1/\sqrt{3})$, we 
have that, for~$n\geq n_1$,
\begin{linenomath}
\begin{align*}
\hP_{n,2n}\Big[\hat{m}(\hat{g}_n) \geq 0\Big] \geq \hP_{n,2n}\Big[\hat{m}(\hat{f}_n) \geq \hat{m}(\hat{f}^2_n)\Big] = \hP_{n,2n}\Bigg[\frac{\hat{m}(\hat{h}_n)}{\sqrt{n}} \geq \delta\frac{\hat{m}(\hat{h}^2_n)}{n}\Bigg].
\end{align*}
\end{linenomath}
On the other hand,
 observe that, under~$\hP_{n,2n}$, 
the random variables $\hat{m}(\hat{h}_n)$ 
and~$\hat{m}(\hat{h}^2_n)$ have the same law 
as~$\hh_n(\hat{X}_1)+\dots+ \hh_n(\hat{X}_{2n})$ and~$\hh^2_n(\hat{X}_1)+\dots+ \hh^2_n(\hat{X}_{2n})$, respectively, where the random variables $\hat{X}_1,\dots,\hat{X}_{2n}$ are i.i.d.~with law~$\hp_n$. Moreover, we have 
that~$|\hh_n(\hat{X}_1)|\leq 1$, $\hp_n$-a.s., $E_{\hp_n}[\hh_n(\hat{X}_1)]=0$ 
and~$\sigma^2:=E_{\hp_n}[\hh^2_n(\hat{X}_1)]\geq 1/2$. If we denote the standard Normal distribution function by~$\Phi$, 
and take~$n_2=n_2(\delta)=4\lceil\frac{1}{(1-\Phi(\delta))^2}\rceil\vee n_1$, then, 
by using the Berry-Esseen theorem (with~$1/2$ as an upper bound 
for the Berry-Esseen constant, see for example~\cite{Tyurin}), 
we obtain that, for~$n\geq n_2$,
\begin{linenomath}
\begin{align*}
\hP_{n,2n}\Bigg[\frac{\hat{m}(\hat{h}_n)}{\sqrt{2n}} \geq \delta\frac{\hat{m}(\hat{h}^2_n)}{n\sqrt{2}}\Bigg] \geq \hP_{n,2n}\Bigg[\frac{\hat{m}(\hat{h}_n)}{\sigma\sqrt{2n}} \geq \frac{\delta}{\sigma\sqrt{2}}\Bigg] \geq c_1,
\end{align*}
\end{linenomath}
where
\begin{linenomath}
\begin{align*}
c_1=c_1(\delta) = \frac{1}{2}\Big(1-\Phi(\delta)\Big).
\end{align*}
\end{linenomath}
Then observe that the above implies that
\begin{linenomath}
\begin{align}
\|\hP_{n,2n}-\hat{\Q}_{n,2n}\|_{\text{TV}} \leq 1-c_1,
\label{TVD_2n}
\end{align}
\end{linenomath}
for all~$n\geq n_2$.

Now, denote by~$\hat{\mu}_{n,k}$ the maximal coupling of~$\hP_{n,k}$ and~$\hat{\Q}_{n,k}$, 
and by $(\hat{m}_1, \hat{m}_2)$ the elements of 
$\hat{\M}_{k}\times \hat{\M}_{k}$. 
Let~$\mathsf{K}$ be the coupling event (that is, $\mathsf{K}=\{(\hat{m}_1, \hat{m}_2): \hat{m}_1 = \hat{m}_2\}$), $\mathsf{K}_1$
 the coupling event of~$(\hat{m}_i(\{\0_1\}), \hat{m}_i(\{\0_2\}))$, 
for $i=1,2$, and $\mathsf{K}_2$ the coupling event of~$(\hat{m}_i(A))_{A\in\mathcal{T}}$, for $i=1,2$, 
and also observe that~$\mathsf{K}=\mathsf{K}_1\cap\mathsf{K}_2$.
 Thus, we deduce that, for all~$n,k\in\IN$,
 \begin{linenomath}
\begin{align*}
\|\hP_{n,k}-\hat{\Q}_{n,k}\|_{\text{TV}} &= 1- \hat{\mu}_{n,k}[\mathsf{K}] \\
&= 1- \sum_{\ell=0}^{k}\hat{\mu}_{n,k}[\mathsf{K}_1 \cap \mathsf{K}_2, \hat{m}_1(\{\0_1,\0_2\})= \hat{m}_2(\{\0_1,\0_2\})=\ell ] \\
&\geq \sum_{\ell=0}^{k} \Big(1-\hat{\mu}_{n,k}[\mathsf{K}_2 \mid \hat{m}_1(\{\0_1,\0_2\})= \hat{m}_2(\{\0_1, \0_2\})=\ell] \Big) \mathfrak{p}_{\ell}^{k} \\
&= \sum_{\ell=0}^{k} \hat{\mu}_{n,k}[\mathsf{K}_2^c \mid \hat{m}_1(\{\0_1,\0_2\})= \hat{m}_2(\{\0_1,\0_2\})=\ell] \mathfrak{p}_{\ell}^{k} \\
&\geq \sum_{\ell=0}^{k} \mathfrak{p}_{\ell}^{k} \|\P_{n,\ell}-\Q_{n,\ell}\|_{\text{TV}},
\end{align*}
\end{linenomath}
where~$\mathfrak{p}_{\ell}^{k}$ is the probability mass function of a Binomial($k,1/2$) random variable at~$\ell$ and $\P_{n,\ell}$ (respectively,~$\Q_{n,\ell}$) is a binomial process with parameters $(\p_n,\ell)$ (respectively,~$(\q_n,\ell)$). 

Using~\eqref{TVD_2n} and the fact 
that~$\|\hP_{n,k}-\hat{\Q}_{n,k}\|_{\text{TV}}$ is non decreasing 
in~$k$ (this follows from the fact 
that~$\|\hP_{n,k}-\hat{\Q}_{n,k}\|_{\text{TV}} 
= \frac{1}{2}E[|1-\mathcal{L}_k(X_1, \dots, X_k)|]$, 
where~$(X_i)_{i\geq 1}$ are i.i.d random variables with 
law~$\p_n$ and~$\mathcal{L}_k(X_1, \dots, X_k)
:=\Pi_{i=1}^k\frac{\mathrm{d}\q_n}{\mathrm{d}\p_n}(X_i)$ 
is a martingale under the canonical filtration), we obtain that,
 for all~$n\geq n_2$ and~$i \leq n$,
 \begin{linenomath}
\begin{align}
\sum_{k=0}^{2i} \mathfrak{p}_{k}^{2i} \|\P_{n,k}-\Q_{n,k}\|_{\text{TV}} 
\leq 1 - c_1.
\label{bound_sum}
\end{align}
\end{linenomath}

Using again the Berry-Esseen theorem 
(once again with~$1/2$ as an upper bound for the Berry-Esseen constant), 
we can deduce that there exist $n_3=n_3(\delta)
=\lceil(\frac{12}{1-\Phi(\delta)})^2\rceil$ and $c_2=c_2(\delta)=-\frac{1}{2}\Phi^{-1}(\frac{1-\Phi(\delta)}{24})\geq1$, 
such that for all $i\geq n_3$, we have
\begin{linenomath}
\begin{align*}
\sum_{k\in \mathcal{J}_i(\delta)} \mathfrak{p}_{k}^{2i} 
\geq 1 - \frac{c_1}{3},
\end{align*}  
\end{linenomath}
where~$\mathcal{J}_i(\delta) := [i-c_2\sqrt{i}, i+c_2\sqrt{i}]$. 
On the other hand, if $i\geq n_2$, by~\eqref{bound_sum} it follows that
\begin{linenomath}
\begin{align*}
\sum_{k\in \mathcal{J}_i(\delta)} 
\mathfrak{p}_{k}^{2i} \|\P_{n,k}-\Q_{n,k}\|_{\text{TV}} 
\leq 1 - c_1(\delta),
\end{align*}
\end{linenomath}
so that, if $i\geq n_2\vee n_3$, there exists~$i_0\in\mathcal{J}_i(\delta)$ such that
\begin{linenomath}
\begin{align*}
\|\P_{n,i_0}-\Q_{n,i_0}\|_{\text{TV}} 
\leq \frac{1 - c_1}{1 - \frac{c_1}{3}} \leq 1 - \frac{2}{3}c_1.
\end{align*}
\end{linenomath}

To conclude the proof of the lemma, observe that for any~$n$ large enough, there exists~$i\geq n_2\vee n_3\vee c_2^2$ such that~$n-(i+\lfloor c_2\sqrt{i} \rfloor)\in[0,3]$, 
the above argument allows to obtain~$i_0\in\IN$ such that~$n-5c_2\sqrt{n} \leq i_0 \leq n$ and
\begin{linenomath}
\begin{align}
\|\P_{n,i_0}-\Q_{n,i_0}\|_{\text{TV}} \leq 1 - \frac{2}{3}c_1.
\label{i0}
\end{align}
\end{linenomath}

Now, if $n-i_0\leq (n_2\vee n_3\vee c_2^2)+5\sqrt{5} c_2^{3/2}n^{1/4}$,
 using (\ref{i0}) and making a point-by-point coupling between 
the $n-i_0$ remaining points of the binomial processes~$\P_{n}$
 and~$\Q_{n}$, we obtain that
 \begin{linenomath}
\begin{align}
\|\P_{n}-\Q_{n}\|_{\text{TV}}  \leq 1-\frac{2}{3}c_1(\delta)\Bigg(1-\frac{[(n_2\vee n_3\vee c_2^2)+5\sqrt{5} c_2^{3/2}n^{1/4}]\delta}{2\sqrt{n}}\Bigg)^+.
\label{Case1}
\end{align}
\end{linenomath}

On the other hand, if 
$n-i_0> (n_2\vee n_3\vee c_2^2)+5\sqrt{5} c_2^{3/2}n^{1/4}$, 
we first consider $j\in \mathbb{N}$ such that 
\begin{linenomath}
\begin{equation*}
n-i_0 - ( j-i_0+ \lfloor c_2\sqrt{j-i_0}\rfloor) \in [0,3].
\end{equation*}
\end{linenomath}
Observe that in this case, $j-i_0> n_2\vee n_3\vee c_2^2$ 
and thus by the former analysis we obtain that there exists $j_0>i_0$ 
such that~$n-5\sqrt{5}c_2^{3/2}n^{1/4} \leq j_0 \leq n$ and
\begin{linenomath}
\begin{align}
\|\P_{n,j_0-i_0}-\Q_{n,j_0-i_0}\|_{\text{TV}} \leq 1 - \frac{2}{3}c_1.
\label{j0}
\end{align}
\end{linenomath}

Using (\ref{i0}), (\ref{j0}) and performing a point-by-point 
coupling between the $n-j_0$ remaining points of the binomial 
processes~$\P_{n}$ and~$\Q_{n}$, we obtain that
\begin{linenomath}
\begin{align}
\|\P_{n}-\Q_{n}\|_{\text{TV}}  \leq 1-\Big(\frac{2}{3}c_1\Big)^2\Bigg(1-\frac{5\sqrt{5} c_2^{3/2}n^{1/4}\delta}{2\sqrt{n}}\Bigg)^+.
\label{Case2}
\end{align}
\end{linenomath}
Finally, using~(\ref{Case1}) and~(\ref{Case2}), 
we obtain that there exists $n_4=n_4(\delta)$ such that, 
for $n\geq n_4$ we have that 
\begin{linenomath}
\begin{align}
\label{OPIO}
\|\P_{n}-\Q_{n}\|_{\text{TV}}  \leq 1-c_3,
\end{align}
\end{linenomath}
where $c_3=c_3(\delta)$ is a positive constant depending on~$\delta$. 

To finish the proof of Lemma~\ref{Propmulti_2}, 
we consider the case $n<n_4$. 
Since $\q_n\ll\p_n$ and 
$|\frac{\mathrm{d}\q_n}{\mathrm{d}\p_n}(x)-1|\leq\delta$
 for all~$x\in \Omega$, we first observe that 
 \begin{linenomath}
\begin{equation*}
\max_{1\leq n< n_4}\|\p_n-\q_n\|_{\text{TV}}\leq 1-\frac{1}{1+\delta}.
\end{equation*}
\end{linenomath}
Then, from this last fact we obtain that
\begin{linenomath}
\begin{equation*}
\max_{1\leq n< n_4} \|\P_{n}-\Q_{n}\|_{\text{TV}}
\leq 1-\Big(\frac{1}{1+\delta}\Big)^{n_4}.
\end{equation*}	
\end{linenomath}
Together with~(\ref{OPIO}), this concludes the proof of
 Lemma~\ref{Propmulti_2}.
\end{proof}

We now prove Theorem~\ref{Thm2}.
In this last part, we consider the following decomposition of the transition density~$p$,
 \begin{linenomath}
	\begin{align*}
	p(x,\cdot) = q + (1-q)\mu'(x,\cdot), \text{ for all } x\in\Sigma,
	\end{align*}
\end{linenomath}
where~$q:=1-\varepsilon$, for~$\varepsilon\in(0,1)$, and~$\mu'(x,\cdot)=\frac{p(x,\cdot) - q}{1-q}\geq 0$ is a probability density with respect to~$\Pi$, since~$\alpha\geq 1-\varepsilon$. Observe that we can construct the same coupling of Section~\ref{coupling}, but now using Bernoulli random variables with parameter~$q$ instead of $\frac{1}{2}$, replacing the event~$\C$ and the quantity~$F$ respectively by
\begin{linenomath}
	\begin{align*}
	\tilde{\mathsf{C}}:=\Big\{|\mH|> \frac{q^2}{6}n\Big\}
	\end{align*}
\end{linenomath}
(recall that~$\mH$ is defined in~\eqref{set_H}) and
\begin{linenomath}
	\begin{align*}
	\tilde{F}:= C_{11} \displaystyle\sqrt{1 + \ln(\phi 2^{\beta}) + \frac{\beta}{\gamma}\ln\Big(\frac{\kappa\vee (2\eps)}{\eps}\Big)},
	\end{align*}
\end{linenomath}
where~$C_{11}=C_{11}(\eps)$, maintaining the same notations for the other quantities defined in that section. Following the same steps as in the proof of Proposition~\ref{Prop_H_comp}, we obtain

{\prop \label{Prop_H_comp_2} There exist $n_5=n_5(\eps)\in \IN$, such that, for all $n\geq n_5$,
	it holds that
	\begin{linenomath}
		\begin{align*}
		\IP\big[\tilde{\mathsf{C}}^c\big] \leq \frac{1}{2}.
		\end{align*}
	\end{linenomath}
}

Then, using Proposition~\ref{Prop_H_comp_2} and the same reasoning as in the proof of Proposition~\ref{goodenv}, we obtain

{\prop  \label{goodenv_2} There exist a positive~$C_{12}=C_{12}(\eps)$ and $n_6=n_6(\eps)\in\IN$,
	such that, for all~$n\geq n_6$ and $i\in\IN$, it holds that
	\begin{linenomath}
		\begin{align*}
		\IP\big[\A_i^c\;\big|\; \tilde{\mathsf{C}} ~\big] &\leq  \frac{C_{12}}{(1+i)^3}.
		\end{align*}
	\end{linenomath} 
}

Then, we take $i=i_1:=\lceil (2C_{12})^{1/3}\rceil-1$ 
in Proposition~\ref{goodenv_2}, so that 
\begin{linenomath}
\begin{equation*}
\IP\big[\A_{i_1} \;\big|\; \tilde{\mathsf{C}} ~\big]
\geq \frac{1}{2},
\end{equation*}
\end{linenomath}
for all~$n\geq n_6$. Then, using Proposition~\ref{Prop_H_comp_2},
 we obtain that $\IP[\A_{i_1}]\geq 1/4$, for $n\geq n_5\vee n_6$. 
Now, observe that, by Lemma~\ref{Propmulti_2}, for all
 $n\geq n_7 = n_7(\beta,\varphi, \kappa, \gamma, \eps) 
:=\max\{n_5, n_6, [\tilde{F}(1+i_1)]^2\}$ we obtain that
\begin{equation}
\label{TIRIO}
\text{d}_{\text{TV}}(L_n^X , L_n^Y)\leq 1-\IP[\Upsilon]
\leq 1-\IP[\Upsilon,\A_{i_1}]\leq 1-\frac{1}{4}C_{10}((1+i_1)\tilde{F}).
\end{equation}
Then, to complete the proof we just observe that
\begin{equation}
\label{TYREL}
\max_{1\leq n< n_7}\text{d}_{\text{TV}}(L_n^X , L_n^Y)
\leq 1-\Big(\frac{1}{1+\eps}\Big)^{n_7}.
\end{equation}
Finally, using~(\ref{TIRIO}), (\ref{TYREL}) and observing that (if needed) $K'$ can always be modified to be decreasing in $\eps$, we conclude
the proof of Theorem~\ref{Thm2}.
\qed

\section*{Acknowledgements}
Diego F.~de Bernardini thanks S\~ao Paulo Research Foundation (FAPESP) 
(grant \#2016/13646--4) and Fundo de Apoio ao Ensino, 
\`a Pesquisa e \`a Extensão (FAEPEX) (grant \#2866/16). Christophe 
Gallesco was partially supported by CNPq (grant 313496/2014--5).
Serguei Popov was partially supported by CNPq (grant 300886/2008--0).
The three authors were partially supported by FAPESP
(grant \#2017/02022--2).
 The authors are very thankful to Caio Alves and the anonymous referee 
for the careful reading of the first version of this paper and useful 
comments and suggestions.

\end{document}